\documentclass[11pt,a4paper]{amsart}
\IfFileExists{lmodern.sty}{\usepackage[T1]{fontenc}\usepackage{lmodern}}{}
\usepackage[a4paper,margin=26mm]{geometry}
\usepackage[utf8]{inputenc}
\usepackage{amsmath,amssymb,amsthm}
\usepackage{enumitem,tabularx,booktabs}
\setlist[enumerate]{label=\textup{(\alph*)},leftmargin=2.3em,labelsep=0.6em,topsep=3pt,itemsep=2pt,parsep=0pt,partopsep=0pt,beginpenalty=10000}
\usepackage{microtype,needspace}
\usepackage{needspace}
\usepackage{hyperref}
\hypersetup{colorlinks=true,linkcolor=black,citecolor=black,urlcolor=black}

\newtheorem{theorem}{Theorem}[section]
\newtheorem{proposition}[theorem]{Proposition}
\newtheorem{lemma}[theorem]{Lemma}
\newtheorem{corollary}[theorem]{Corollary}
\theoremstyle{definition}
\newtheorem{definition}[theorem]{Definition}
\newtheorem{example}[theorem]{Example}
\theoremstyle{remark}
\newtheorem{remark}[theorem]{Remark}
\numberwithin{equation}{section}

\newcommand{\N}{\mathbb{N}}
\newcommand{\Z}{\mathbb{Z}}
\newcommand{\F}{\mathbb{F}}
\newcommand{\Q}{\mathbb{Q}}
\newcommand{\rad}[1]{\sqrt{#1}}
\DeclareMathOperator{\Min}{Min}
\DeclareMathOperator{\Ann}{Ann}
\DeclareMathOperator{\Char}{char}
\DeclareMathOperator{\Nil}{Nil}

\begin{document}
\raggedbottom

\title[$(m,n)$-mpdf ideals]{$(m,n)$-mixed-power difference factor absorbing ideals of commutative rings}

\author{Omar H. Taha}
\address{Department of Mathematics, College of Computer Science and Mathematics, Tikrit University, Tikrit, Iraq.}
\email{omar.h.tahamm2314@st.tu.edu.iq}

\author{Akram Salim}
\address{Department of Mathematics, College of Computer Science and Mathematics, Tikrit University, Tikrit, Iraq.}
\email{akr\_tel@tu.edu.iq}

\author{Hwankoo Kim}
\address{Division of Computer Engineering, Hoseo University, Republic of Korea.}
\email{hkkim@hoseo.edu}


\begin{abstract}
Let $n\ge2$ and $k\ge1$, and put $m=kn$. A proper ideal $I$ of a commutative ring $R$ is $(m,n)$-mixed-power difference factor absorbing if for all $a,b \in R$ with $a^{kn}-b^n\in I$ implies that $a^k-b$ or $F_{k,n} \in I$, where $F_{k,n}(a,b)=\sum_{i=0}^{n-1}a^{k(n-1-i)}b^i$.
We study how this condition behaves under ring constructions. 
Our principal result gives an exact criterion for arbitrary direct products. 
We show that the requirement that every mixed-power identity select the linear factor is equivalent to injectivity of the $n$th-power map, independently of $k$.
This gives an explicit product classification
and a connectedness obstruction for nonradical mpdf quotients on which $n$ is a non-zero-divisor.
For even $n$, that regularity hypothesis forces every mpdf ideal to be primary. 
For a $P$-primary ideal with $2\notin P$, the strict quartic condition and the fourth-power
pair-product condition are both equivalent to $P^3\subseteq I$.
The proof gives a cube-zero quotient nilradical whenever the pair-product condition holds and $2$ is regular on $R/I$. 
We also establish localization and chain properties, extend the square-case characterization to every multiplier $k$, and relate prime-power exponents in characteristic $p$ to closed ideals. 
Examples show genuine dependence on $k$ for nonradical ideals and failure of polynomial ascent and local detection. Pullbacks, idealizations and amalgamated algebras provide applications of the transfer and product criteria.
\end{abstract}

\maketitle

\section{Introduction}\label{sec:intro}

Throughout, all rings are commutative with identity. An ideal $I$ of a ring $R$ is prime if $I\ne R$ and $xy\in I$ implies $x\in I$ or $y\in I$. Many generalizations of prime ideals keep this requirement only for products of a special shape. One such class was introduced by Anderson, Badawi and Coykendall~\cite{ABC}: a proper ideal $I$ is square-difference factor absorbing, or sdf-absorbing, if $a^2-b^2=(a-b)(a+b)\in I$ implies $a-b\in I$ or $a+b\in I$ for all nonzero $a,b\in R$. Farshadifar~\cite{Far25} replaced squares by cubes, using $a^3-b^3=(a-b)(a^2+ab+b^2)$, and Ye\c{s}ilot~\cite{Yes26} treated an arbitrary prime exponent $p$. Hizem~\cite{Hizem} described sdf-absorbing ideals through their minimal primes, and submodule, primary and semiprimary versions of the square and cube conditions have also been considered~\cite{KC25,KCT,Dr26}.

In this paper we look at differences of powers with different exponents. For integers $n\ge2$ and $k\ge1$ put $m=kn$. Then
\[
a^m-b^n=(a^k-b)\,F_{k,n}(a,b),\qquad F_{k,n}(a,b)=\sum_{i=0}^{n-1}a^{k(n-1-i)}b^i.
\]
We call a proper ideal $I$ of $R$ an \emph{$(m,n)$-mixed-power difference factor absorbing ideal}, or \emph{$(m,n)$-mpdf ideal}, if $a^m-b^n\in I$ implies $a^k-b\in I$ or $F_{k,n}(a,b)\in I$ for all $a,b\in R$. For $k=1$ this is an equal-power condition, and $(2,2)$, $(3,3)$ and $(p,p)$ correspond to the sdf, cdf and prime-exponent conditions above. Since $F_{k,n}(a,b)=F_{1,n}(a^k,b)$, the element $a$ enters only through $a^k$, and the condition does depend on $k$: the ideal $(16)\subset\Z$ is $(6,3)$-mpdf but not $(3,3)$-mpdf.

We allow $a$ and $b$ to be zero, as in the published cubic definition \cite[Definition 2.1]{Far25} and the prime-exponent definition of \cite{Yes26}. Thus our $(3,3)$ and $(p,p)$ cases are those notions. 
The sdf condition of~\cite{ABC} excludes zero. The two square conditions agree for nonzero ideals but differ for the zero ideal of $\Z_4$ (Proposition~\ref{prop:sdfcdf}); we retain this distinction throughout.
Ye\c{s}ilot, \c{S}im\c{s}ek and \"Ozel Ay~\cite{YSO} study \emph{cdf-primary} ideals, requiring $a-b\in I$ or $a^2+ab+b^2\in\sqrt I$, and give exact finite-product criteria \cite[Propositions 14 and 15]{YSO}. 
Draoui's semiprimary version places both conclusions in $\sqrt I$ \cite[Definition 2.1]{Dr26}. 
These cubic conditions agree on radical ideals but differ in general.
Their product arguments already identify the conflict between factors chosen in different coordinates; Section~\ref{sec:products} explains the precise relation to our criterion, including nonradical factors.

For fourth powers, Tekir, Anebri, Hasi\c{c} and Arabac{\i} \cite{TAHA} select a product of two of the three factors $(a-b),(a+b),(a^2+b^2)$. 
A proper ideal $I$ of $R$ is \emph{$4$-power-difference-factor absorbing}, or \emph{$4$-pdf-absorbing}, if, for all $a,b\in R$, with $a^4-b^4\in I$ implies 
\[
 a^2-b^2\in I,\quad\text{or}\quad
 (a^2+b^2)(a-b)\in I,\quad\text{or}\quad
 (a^2+b^2)(a+b)\in I.
\].
Their $4$-pdf-absorbing condition is weaker than our $(4,4)$-mpdf condition, which selects $a-b$ or the full
cubic cofactor. Both definitions quantify over all pairs. We give the precise comparison in Definition of $4$-pdf and Proposition~\ref{quartic:implication}.
The primary criterion of \cite[Theorem~4.6]{TAHA} assumes $2,3\notin P$ and characterizes $4$-pdf-absorbing $P$-primary ideals by $P^3\subseteq I$. 
We retain only $2\notin P$ and show that the same criterion also characterizes the strict quartic condition.

Another related notion is due to Anderson and Badawi~\cite{AB17}: for positive integers $u,v$, an ideal $I$ is \emph{$(u,v)$-closed} if $c^u\in I$ implies $c^v\in I$ for every $c\in R$. Every $(m,n)$-mpdf ideal is $(n,n-1)$-closed, and in prime characteristic this necessary condition is often sufficient. The notation $(m,n)$ is also used for absorbing ideals~\cite{Zarei,BadawiKhalfiMahdou}, but those conditions are different from ours.

The structural question is whether the two factors can be selected consistently across a quotient. In a direct product, a power identity holds in each coordinate, but different coordinates can force opposite factors. Nilpotents create a further obstruction: a nonzero element can have zero $n$-th power even when a power map is injective on the residue field. 
These observations explain why quotient transfer alone does not decide the condition, and why the multiplier $k$ deserves separate study.

The main result extends the finite cubic product mechanism of \cite[Propositions 14 and 15]{YSO} to arbitrary exponents $n$, mixed multipliers $k$, and arbitrary index sets (Theorem~\ref{thm:product}). Its finite $(3,3)$ specialization is closely related to those criteria even for nonradical factors.
Proposition~\ref{prop:linear-injectivity} identifies the universal linear-factor alternative with injectivity of the $n$-th-power map. Together they give the explicit classification in Corollary~\ref{cor:product-explicit} and the connectedness criterion in Corollary~\ref{cor:connected-mpdf}.
For even $n$, regularity of $n$ modulo an mpdf ideal forces the ideal to be primary (Theorem~\ref{thm:even-primary}); examples show why neither evenness nor the distinction from primeness can be omitted.
Theorem~\ref{quartic:primary-equivalence} gives the quartic primary comparison. Its proof uses the individual nilpotency bound from \cite[Proposition~4.1]{TAHA} and an additional pair of substitutions
to obtain $(\sqrt I)^3\subseteq I$ with $2$ regular modulo $I$, without a hypothesis on $3$. 
The square case (Theorem~\ref{thm:n2}) is independent of $k$, whereas Example~\ref{ex:16} exhibits a strict change between $k=1$ and $k=2$.
The square results extend, and explicitly retain attribution to, the work of Anderson, Badawi and Coykendall~\cite{ABC} and Hizem~\cite{Hizem}. 
The quotient and canonical-construction statements are foundational transfer results; their role is to support the structural criteria.

Section~\ref{sec:def} gives the definition, quotient formulation and examples. Section~\ref{sec:structure} treats chains, the square case, primary ideals for even exponents, the quartic primary comparison, and prime-power exponents in positive characteristic. Section~\ref{sec:hom} contains the localization, product and connectedness results. 
Section~\ref{sec:classical} compares the condition with familiar classes of ideals, 
and Section~\ref{sec:poly} records the limits of polynomial and power series transfer.
Finally, Section~\ref{sec:constructions} applies the results to pullbacks, idealizations~\cite{AndersonWinders} and amalgamated algebras~\cite{DAnna}. 
The companion paper~\cite{MPDFII} develops the general spectral and arithmetic classifications, the radical quartic comparison, and the module and extension data needed for square-zero quotients.

\section{Definition, reformulations and first examples}\label{sec:def}

Throughout, rings are commutative with $1\ne0$, and homomorphisms preserve the identity; a zero ring may occur only as an auxiliary quotient in a pullback. Unless specified otherwise, ideals are proper, $n\ge2$, $k\ge1$ and $m=kn$. We write $\Z_N=\Z/N\Z$, $\rad I=\{x:x^t\in I\text{ for some }t\ge1\}$, $\Min_R(I)$ for the set of primes of $R$ minimal over $I$, $R^\times$ for the group of units, $\mu_n(K)=\{\zeta\in K:\zeta^n=1\}$ for a field $K$, and $v_p(n)$ for the exponent of a prime $p$ in an integer $n\ne0$. For an integer $n$ and an ideal $J$, the statement $n\in J$ means $n\cdot1\in J$. We use the conventions $x^0=1$ for every element $x$, including $x=0$, and $F_{k,1}(a,b)=1$. Quotient rings, radicals, localization and the Chinese remainder theorem (CRT) are used as in~\cite{AM}.
We first record elementary identities for $F_{k,n}$. They follow by
expanding the finite geometric sum and grouping its terms.

\begin{lemma}\label{lem:poly}
For all $a,b,t$ in a commutative ring and all integers $q,r\ge1$:
\begin{enumerate}
\item $a^m-b^n=(a^k-b)F_{k,n}(a,b)$ and $F_{k,n}(a,b)=F_{1,n}(a^k,b)$.
\item $F_{k,n}(ta,t^kb)=t^{k(n-1)}F_{k,n}(a,b)$.
\item $F_{k,n}(a,0)=a^{k(n-1)}$, $F_{k,n}(0,b)=b^{n-1}$ and $F_{k,n}(a,a^k)=na^{k(n-1)}$.
\item $F_{1,qr}(a,b)=F_{1,q}(a,b)\,F_{1,r}(a^q,b^q)$.
\end{enumerate}
\end{lemma}
The next identities follow from the Frobenius homomorphism in
characteristic $p$. We express them over $\mathbb Z[X,Y]$ for use
modulo ideals containing $p$.

\begin{lemma}\label{lem:charp}
Let $p$ be prime and $n=p^e$ with $e\ge1$. Then there are $G,H\in\Z[X,Y]$ such that
\[
F_{1,n}(X,Y)=(X-Y)^{n-1}+pG\qquad\text{and}\qquad(X-Y)^n=X^n-Y^n+pH.
\]
Consequently, in every commutative ring $R$ and for all $a,b\in R$, $F_{k,n}(a,b)-(a^k-b)^{n-1}\in pR$ and $(a^k-b)^n-(a^m-b^n)\in pR$.
\end{lemma}

\begin{definition}\label{def:mpdf}
Let $m,n\in\N$, $n\ge2$ and $n\mid m$, and write $m=kn$ with $k\ge1$. A proper ideal $I$ of $R$ is called an \emph{$(m,n)$-mixed-power difference factor absorbing ideal}, briefly an \emph{$(m,n)$-mpdf ideal}, if, for all $a,b\in R$, $a^m-b^n\in I$ implies $a^k-b\in I$ or $F_{k,n}(a,b)\in I$. Here
\[
a^m-b^n=(a^k-b)F_{k,n}(a,b)\quad\text{and}\quad F_{k,n}(a,b)=\sum_{i=0}^{n-1}a^{k(n-1-i)}b^i.
\]
\end{definition}

\begin{remark}\label{rem:def}
For $n=1$ the cofactor is $1$, so $a^m-b\in I$ already gives the first alternative and every proper ideal would satisfy the condition. The restriction $n\ge2$ excludes this vacuous case.
\end{remark}

\begin{proposition}\label{prop:reform}
Let $I$ be a proper ideal of $R$ and $S=R/I$. The following statements are equivalent.
\begin{enumerate}[label=\textup{(\roman*)}]
\item $I$ is an $(m,n)$-mpdf ideal of $R$.
\item For all $\alpha,\beta\in S$ with $\alpha^m=\beta^n$ we have $\alpha^k=\beta$ or $F_{k,n}(\alpha,\beta)=0$; that is, the zero ideal of $S$ is $(m,n)$-mpdf.
\end{enumerate}
\end{proposition}
\begin{proof}
For the quotient map $\pi:R\to R/I$, polynomial evaluation commutes
with $\pi$. Thus $a^{kn}-b^n\in I$ is equivalent to
$\pi(a)^{kn}=\pi(b)^n$, and the two ideal-membership conclusions are
equivalent to the corresponding zero equalities. Every pair in $R/I$
has lifts in $R$, so the implications are equivalent.
\end{proof}

\begin{proposition}\label{prop:prime}
For every admissible pair $(m,n)$, the following statements hold.
\begin{enumerate}
\item Every prime ideal of $R$ is an $(m,n)$-mpdf ideal of $R$.
\item If $x\mapsto x^n$ is injective on $R/I$, then $I$ is $(kn,n)$-mpdf for every $k\ge1$.
\end{enumerate}
\end{proposition}
\begin{proof}
\textup{(a)} Let $P$ be a prime ideal and $a,b\in R$ with $a^m-b^n\in P$. By Lemma~\ref{lem:poly}(a), $(a^k-b)F_{k,n}(a,b)\in P$, so $a^k-b\in P$ or $F_{k,n}(a,b)\in P$.

\noindent\textup{(b)} Let $a,b\in R$ with $a^{kn}-b^n\in I$. In $R/I$ this says $(\bar a^k)^n=\bar b^n$. By injectivity of $z\mapsto z^n$, $\bar a^k=\bar b$, that is, $a^k-b\in I$.
\end{proof}

\begin{proposition}\label{prop:kdiv}
Let $I$ be a $(kn,n)$-mpdf ideal and let $k'$ be a positive multiple of $k$. Then $I$ is a $(k'n,n)$-mpdf ideal. In particular, an $(n,n)$-mpdf ideal is $(m,n)$-mpdf for every multiple $m$ of $n$.
\end{proposition}
\begin{proof}
Let $a,b\in R$ and suppose that $a^{k'n}-b^n\in I$. Write $k'=k\ell$ and put $c=a^\ell$. Then
\[
c^{kn}-b^n=a^{k'n}-b^n\in I.
\]
Since $I$ is $(kn,n)$-mpdf, $c^k-b\in I$ or $F_{k,n}(c,b)\in I$. Here $c^k=a^{k'}$ and
\[
F_{k,n}(c,b)=\sum_{i=0}^{n-1}a^{k'(n-1-i)}b^i=F_{k',n}(a,b),
\]
so $I$ is $(k'n,n)$-mpdf.
\end{proof}

The converse of Proposition~\ref{prop:kdiv} fails (Example~\ref{ex:16}).

\begin{example}\label{ex:16}
The ideal $(16)\subset\Z$ is $(3k,3)$-mpdf if and only if $k\ge2$.

\emph{Case $k\ge2$.} Let $a,b\in\Z$ with $16\mid a^{3k}-b^3$.
\begin{itemize}
\item Suppose $a$ is odd. Then $a^{3k}$ and hence $b^3$ is odd, so $b$ is odd. The group $(\Z_{16})^\times$ has order $8$, which is prime to $3$, so cubing is injective on it. Since $(a^k)^3\equiv b^3\pmod{16}$ with $a^k,b$ odd, we get $a^k\equiv b\pmod{16}$.
\item Suppose $a$ is even. Then $16\mid a^{3k}$ because $3k\ge6$, so $16\mid b^3$. If $b\equiv2\pmod4$ then $b^3\equiv8\pmod{16}$, which is impossible; if $b$ is odd then $b^3$ is odd. So $4\mid b$. Also $4\mid a^k$ since $k\ge2$. Hence $a^{2k}$, $a^kb$ and $b^2$ are all divisible by $16$, and $16\mid a^{2k}+a^kb+b^2=F_{k,3}(a,b)$.
\end{itemize}

\begin{proposition}
Let $I$ be a proper ideal and suppose that every element of $R/I$ is a $k$th power in $R/I$. Then $I$ is $(kn,n)$-mpdf if and only if it is $(n,n)$-mpdf.
\end{proposition}
\begin{proof}
If $I$ is $(n,n)$-mpdf, Proposition~\ref{prop:kdiv} gives the $(kn,n)$-mpdf condition.
Conversely, suppose $I$ is $(kn,n)$-mpdf. Let $a,b\in R$ with $a^n-b^n\in I$. By hypothesis there is $\overline{c}\in R/I$ such that $\overline{a}=\overline{c}^k=\overline{c^k}$ thus, $c^k-a\in I$. Thus
\[
c^{kn}-b^n\equiv a^n-b^n\equiv0\pmod I.
\]
Hence $c^k-b\in I$ or $F_{k,n}(c,b)\in I$. Since
\[
c^k-b\equiv a-b\pmod I,\qquad
F_{k,n}(c,b)\equiv F_{1,n}(a,b)\pmod I,
\]
$a-b\in I$ or $F_{1,n}(a,b)\in I$. Therefore $I$ is $(n,n)$-mpdf.
\end{proof}

In particular, the conclusion holds when $R/I$ is a finite field with
$q$ elements and $\gcd(k,q-1)=1$, since raising to the $k$th power
then permutes the nonzero elements. The next statement extends the
cubic power-closure result of~\cite[Proposition 2.4]{Far25} and is
the source of most restrictions on mpdf ideals.

\begin{proposition}\label{prop:special}
Let $I$ be an $(m,n)$-mpdf ideal of $R$. If $c\in R$ and $c^{n}\in I$, then $c^{n-1}\in I$. 
\end{proposition}
\begin{proof}
Let $c\in R$ with $c^n\in I$. Take $a=0$ and $b=c$. Then $a^m-b^n=-c^n\in I$, so by the mpdf property $a^k-b=-c\in I$ or $F_{k,n}(0,c)\in I$. By Lemma~\ref{lem:poly}(c), $F_{k,n}(0,c)=c^{n-1}$. If $c^{n-1}\in I$ we are done. If $-c\in I$, then $c\in I$ and $c^{n-1}=c\cdot c^{n-2}\in I$, because $n-1\ge1$. In either case $c^{n-1}\in I$.
\end{proof}

\begin{corollary}\label{cor:closedradical}
Every $(m,n)$-mpdf ideal is $(n,n-1)$-closed. Moreover, every $(2k,2)$-mpdf ideal is radical.
\end{corollary}
\begin{proof}
The first assertion is Proposition~\ref{prop:special}. Let $I$ be $(2k,2)$-mpdf, and let $c\in\rad I$, say $c^t\in I$. Choose $j\ge0$ with $2^j\ge t$. Then $c^{2^j}=c^t\,c^{2^j-t}\in I$. If $j\ge1$, Proposition~\ref{prop:special} applied to $c^{2^{j-1}}$ gives $c^{2^{j-1}}\in I$, and repeating this step $j$ times we obtain
\[
c^{2^j}\in I\ \Rightarrow\ c^{2^{j-1}}\in I\ \Rightarrow\ \cdots\ \Rightarrow\ c^2\in I\ \Rightarrow\ c\in I.
\]
Hence $\rad I\subseteq I$, and $I=\rad I$.
\end{proof}

\Needspace{4\baselineskip}\subsection{Special cases and examples}\label{sec:examples}

For small parameters the polynomial $F_{k,n}$ is explicit:
\begin{equation*}
\begin{aligned}
F_{1,2}(a,b)&=a+b, & F_{2,2}(a,b)&=a^2+b, & F_{3,2}(a,b)&=a^3+b,\\
F_{1,3}(a,b)&=a^2+ab+b^2, & F_{2,3}(a,b)&=a^4+a^2b+b^2, & F_{1,4}(a,b)&=a^3+a^2b+ab^2+b^3.
\end{aligned}
\end{equation*}
Thus an ideal is $(2,2)$-mpdf if $a^2-b^2\in I$ implies $a-b\in I$ or $a+b\in I$; it is $(4,2)$-mpdf if $a^4-b^2\in I$ implies $a^2-b\in I$ or $a^2+b\in I$; and it is $(3,3)$-mpdf if $a^3-b^3\in I$ implies $a-b\in I$ or $a^2+ab+b^2\in I$.

\begin{proposition}\label{prop:sdfcdf}
Let $I$ be a proper ideal of $R$.
\begin{enumerate}
\item If $I\neq 0$, then $I$ is $(2,2)$-mpdf if and only if $I$ is sdf-absorbing.
\item If the zero ideal is $(2,2)$-mpdf, then it is sdf-absorbing.
\item $I$ is $(3,3)$-mpdf if and only if $I$ is cdf-absorbing, and, for a prime $p$, $I$ is $(p,p)$-mpdf if and only if $I$ satisfies the prime-exponent condition of~\cite{Yes26}; both conditions are quantified over all $a,b\in R$.
\end{enumerate}
\end{proposition}
\begin{proof}
\textup{(a)} Let $I\ne0$. If $I$ is $(2,2)$-mpdf, then in particular the implication holds for $a,b\ne0$, so $I$ is sdf-absorbing. Conversely, assume that $I$ is sdf-absorbing, and fix $0\ne i\in I$. Let $a,b\in R$ with $a^2-b^2\in I$; we must show $a-b\in I$ or $a+b\in I$. There are four cases.
\begin{enumerate}
\item[Case (1)] If $a\ne0\ne b$, this is the sdf condition.
\item[Case (2)] If $a=0\ne b$, then $b^2=-(a^2-b^2)\in I$. Since $i^2\in I$, also $b^2-i^2\in I$, with $b\ne0\ne i$. The sdf condition gives $b-i\in I$ or $b+i\in I$, and in both cases $b\in I$ because $i\in I$. Thus $a-b=-b\in I$.
\item[Case (3)] If $b=0\ne a$, then $a^2\in I$, and the same argument with $a$ in place of $b$ gives $a\in I$, so $a-b=a\in I$.
\item[Case (4)] If $a=b=0$, then $a-b=0\in I$.
\end{enumerate}

\noindent\textup{(b)} Suppose the zero ideal is $(2,2)$-mpdf. If $a,b\ne0$ and $a^2-b^2=0$, then $a-b=0$ or $a+b=0$, so the zero ideal is sdf-absorbing. 

\noindent\textup{(c)} For $k=1$ and $n=3$, respectively $n=p$,
Definition~\ref{def:mpdf} is the condition of
\cite[Definition 2.1]{Far25}, respectively~\cite{Yes26}, with the
same quantifiers.
\end{proof}
The converse of Proposition~\ref{prop:sdfcdf} $(2)$ fails (Example~\ref{ex:22}).

\begin{example}\label{ex:22}
The zero ideal of $\mathbb Z_4$ is sdf-absorbing, but it is not $(2,2)$-mpdf.
Clear $I=(0)$ is sdf-absorbing in $\mathbb Z_4$, see  \cite{ABC}.
We now show that $(0)$ is not $(2,2)$-mpdf, let $a=2$ and $b=0$. Then
\[
a^2-b^2=4-0=0\in(0),
\qquad \text{but} \qquad
a-b=2\ne0,\qquad F_{1,2}(a,b)=a+b=2\ne0 .
\]
Hence, $(0)$ is not $(2,2)$-mpdf.
\end{example}

These are the three pair products of the factorization
$a^4-b^4=(a-b)(a+b)(a^2+b^2)$. To distinguish the two notions, we
also call $(4,4)$-mpdf absorption the \emph{strict quartic condition}.
As in Proposition~\ref{prop:reform}, $I$ is $4$-pdf-absorbing if and
only if the zero ideal of $R/I$ is; see \cite[Corollary~3.3]{TAHA}.

\begin{proposition}\label{quartic:implication}
Every $(4,4)$-mpdf ideal is $4$-pdf-absorbing, and every
$2$-absorbing ideal is $4$-pdf-absorbing. Neither converse holds.
\end{proposition}
\begin{proof}
The linear alternative $a-b\in I$ gives $a^2-b^2\in I$, while the
cubic alternative is already one of the $4$-pdf conclusions.
For a $2$-absorbing ideal, apply its defining condition
\cite{Ba07} to the triple $(a-b),(a+b),(a^2+b^2)$; this implication
is also observed immediately after \cite[Definition~2.1]{TAHA}.

The zero ideal of $\Q\times\Q$ is $4$-pdf-absorbing because equal
fourth powers force equal squares in each coordinate. However,
$a=(1,1)$ and $b=(1,-1)$ have equal fourth powers, while
$a-b=(0,2)\ne0$ and $(a+b)(a^2+b^2)=(4,0)\ne0$.
The same pair disproves $(4k,4)$-mpdf absorption for every $k\ge1$.
Similarly, the zero ideal of $\Q^3$ is $4$-pdf-absorbing but is not
$2$-absorbing: the elements $(1,1,0),(1,0,1),(0,1,1)$ have zero
product and three nonzero pair products. Other examples separating
$4$-pdf absorption from $2$-absorption appear in
\cite[Examples~2.2 and~2.6]{TAHA}.
\end{proof}

\begin{example}\label{ex:six}
The ideal $(6)\subset\Z$ is $(2k,2)$-mpdf for every $k\ge1$ and is not prime.

Let $a,b\in\Z$ and suppose that $6\mid a^{2k}-b^2$. Put $x=a^k$, so that $6\mid x^2-b^2=(x-b)(x+b)$.
\begin{itemize}
\item[Case (1)] Modulo $2$: $x^2\equiv b^2$ forces $x\equiv b\pmod 2$, so $x-b$ and $x+b$ are both even.
\item[Case (2)] Modulo $3$: $3$ is prime and divides $(x-b)(x+b)$, so $3\mid x-b$ or $3\mid x+b$.
\end{itemize}
Hence $6\mid x-b=a^k-b$ or $6\mid x+b=F_{k,2}(a,b)$. Finally $2\cdot3=6\in(6)$ while $2,3\notin(6)$, so $(6)$ is not prime. 
\end{example}
Theorem~\ref{thm:n2} gives the same conclusion for every $k$ without computation.
\begin{example}
The ideal $(4)\subset\Z$ is $(3k,3)$-mpdf for every $k\ge1$, but is not radical.

By Proposition~\ref{prop:kdiv} it suffices to treat $k=1$. The cubes of $0,1,2,3$ modulo $4$ are $0,1,0,3$. Let $a,b\in\Z$ with $4\mid a^3-b^3$.
\begin{itemize}
\item[Case (1)] If $a$ is odd, then $a^3\equiv a\pmod 4$, so $b^3\equiv a^3$ is odd, hence $b$ is odd and $b^3\equiv b\pmod4$. Therefore $a\equiv a^3\equiv b^3\equiv b\pmod4$, that is, $4\mid a-b$.
\item[Case (2)] If $a$ is even, then $a^3\equiv0\pmod4$, so $b^3\equiv0\pmod4$ and $b$ is even. Then $a^2$, $ab$ and $b^2$ are each divisible by $4$, so $4\mid a^2+ab+b^2=F_{1,3}(a,b)$.
\end{itemize}
Thus $(4)$ is $(3,3)$-mpdf. It is not radical, since $2^2\in(4)$ and $2\notin(4)$.
\end{example}

\emph{Case $k=1$.} Take $a=2,b=6$. Then $a^3-b^3=8-216=-208=-13\cdot16\in(16)$, while $a-b=-4\notin(16)$ and $a^2+ab+b^2=4+12+36=52\notin(16)$. Hence $(16)$ is not $(3,3)$-mpdf.
\end{example}

\begin{example}\label{ex:notmpdf}
\begin{enumerate}
\item The radical ideal $(15)\subset\Z$ is not $(2,2)$-mpdf: for $a=4,b=1$ we have $a^2-b^2=15\in(15)$, while $a-b=3\notin(15)$ and $a+b=5\notin(15)$.
\item The primary ideal $(9)\subset\Z$ is not $(3,3)$-mpdf: for $a=1,b=-2$ we have $a^3-b^3=1+8=9\in(9)$, while $a-b=3\notin(9)$ and $a^2+ab+b^2=1-2+4=3\notin(9)$.
\item The ideal $(0)\subset\Z\times\Z$ is not $(2k,2)$-mpdf for any $k\ge1$: for $a=(1,1)$ and $b=(1,-1)$ we have $a^{2k}=(1,1)=b^2$, while $a^k-b=(0,2)\ne0$ and $a^k+b=(2,0)\ne0$.
\end{enumerate}
\end{example}

\section{Structural properties}\label{sec:structure}

The necessary power-closure condition controls the square case and
prime-power exponents. We first record order properties of the class,
then identify when the square condition reduces to primeness.

\begin{remark}
In $S=R/I$, the mpdf condition says that there is no pair $a,b$ for
which $a^k-b$ and $F_{k,n}(a,b)$ are both nonzero and their product is
zero. Equivalently, for $a^k\ne b$, membership of $F_{k,n}(a,b)$ in
$\Ann_S(a^k-b)$ forces $F_{k,n}(a,b)=0$.
\end{remark}

\begin{proposition}\label{prop:chain}
Let $\{I_\lambda\}$ be a nonempty family of $(m,n)$-mpdf ideals of $R$ which is totally ordered by inclusion. Then the following ideals are $(m,n)$-mpdf.
\begin{enumerate}
\item $U=\bigcup_\lambda I_\lambda$.
\item $V=\bigcap_\lambda I_\lambda$.
\end{enumerate}
\end{proposition}
\begin{proof}

\noindent\textup{(a)} Eventuality, $V$ is proper ideal of $R$.  Let $a,b\in R$ with $a^m-b^n\in U$. Choose $\lambda$ with $a^m-b^n\in I_\lambda$. Since $I_\lambda$ is $(m.n)$-mpdf, $a^k-b\in I_\lambda\subseteq U$ or $F_{k,n}(a,b)\in I_\lambda\subseteq U$.

\noindent\textup{(b)} Let $a,b\in R$ with $a^m-b^n\in V$, and suppose, for a contradiction, that $a^k-b\notin V$ and $F_{k,n}(a,b)\notin V$. Choose $\lambda$ with $a^k-b\notin I_\lambda$ and $\mu$ with $F_{k,n}(a,b)\notin I_\mu$, and let $I_\nu$ be the smaller of $I_\lambda$ and $I_\mu$. Then neither $a^k-b$ nor $F_{k,n}(a,b)$ lies in $I_\nu$, since $I_\nu$ is contained in both $I_\lambda$ and $I_\mu$. But $a^m-b^n\in V\subseteq I_\nu$, which contradicts the fact that $I_\nu$ is $(m.n)$-mpdf. Hence $a^k-b\in V$ or $F_{k,n}(a,b)\in V$.
\end{proof}

\begin{corollary}
Every $(m,n)$-mpdf ideal contains an $(m,n)$-mpdf ideal which is minimal among $(m,n)$-mpdf ideals. 
\end{corollary}
\begin{proof}
Let $I$ be mpdf and let $\mathcal S$ be the set of $(m,n)$-mpdf ideals contained in $I$, partially ordered by reverse inclusion. $\mathcal S$ is nonempty since $I\in\mathcal S$. If $\mathcal C\subseteq\mathcal S$ is a nonempty chain, then $\bigcap_{J\in\mathcal C}J$ is $(m.n)$-mpdf by Proposition~\ref{prop:chain}(b), and it lies in $\mathcal S$ and bounds $\mathcal C$. By Zorn's lemma, $\mathcal S$ has a maximal element for reverse inclusion, that is, an mpdf ideal contained in $I$ that contains no smaller mpdf ideal. It is minimal among all mpdf ideals, since any mpdf ideal contained in it is also contained in $I$.

\end{proof}

\begin{example}
The minimal mpdf ideals need not be prime: the zero ideal of $\Z_6$ is $(2,2)$-mpdf.
Since $(6)\subset\Z$ is $(2,2)$-mpdf by Example~\ref{ex:six}, so the zero ideal of $\Z/(6)=\Z_6$ is $(2,2)$-mpdf by Proposition~\ref{prop:reform}. It is the smallest ideal of $\Z_6$, hence minimal among mpdf ideals, and it is not prime.
\end{example}

\Needspace{4\baselineskip}\subsection{The square case and primeness}

\begin{lemma}\label{lem:n2}
Let $I$ be a $(2k,2)$-mpdf ideal of $R$ and let $x,y\in R$ with $xy\in I$. Then $2x\in I$ or $2y\in I$.
\end{lemma}
\begin{proof}
Every mixed term in the binomial expansion of $(x+y)^k$ contains the factor $xy\in I$, hence $(x+y)^k\equiv x^k+y^k\pmod I$. Put $a=x+y$ and $b=y^k-x^k$. Then
\begin{align*}
a^k-b&\equiv(x^k+y^k)-(y^k-x^k)=2x^k\pmod I,\\
a^k+b&\equiv(x^k+y^k)+(y^k-x^k)=2y^k\pmod I,
\end{align*}
and therefore
\[
a^{2k}-b^2=(a^k-b)(a^k+b)\equiv4x^ky^k=4(xy)^k\equiv0\pmod I.
\]
Since $F_{k,2}(a,b)=a^k+b$, the mpdf property gives $a^k-b\in I$ or $a^k+b\in I$, that is, $2x^k\in I$ or $2y^k\in I$. If $2x^k\in I$, then $(2x)^k=2^{k-1}\cdot2x^k\in I$, and $2x\in I$ because $I$ is radical by Corollary~\ref{cor:closedradical}. In the same way $2y^k\in I$ gives $2y\in I$.
\end{proof}

\begin{theorem}\label{thm:prime}
Let $I$ be a $(2k,2)$-mpdf ideal of $R$.
\begin{enumerate}
\item If $2$ is a non-zero-divisor on $R/I$ \textup{(}in particular, if $2$ is a unit modulo $I$\textup{)}, then $I$ is prime.
\item If $I$ is primary, then $I$ is prime.
\end{enumerate}
\end{theorem}
\begin{proof}
\textup{(a)} Let $x,y\in R$ with $xy\in I$. By Lemma~\ref{lem:n2}, $2x\in I$ or $2y\in I$, say $2x\in I$. Then $(2+I)(x+I)=0$ in $R/I$, and since $2+I$ is a non-zero-divisor, $x\in I$. Similarly, $2y\in I$ gives $y\in I$.

\noindent\textup{(b)} Let $x,y\in R$ with $xy\in I$ and $x\notin I$. Because $I$ is primary, $y^t\in I$ for some $t\ge1$, so $y\in\rad I$. By Corollary~\ref{cor:closedradical}, $\rad I=I$, so $y\in I$.
\end{proof}

\begin{remark}
The assumptions in Theorem~\ref{thm:prime} matter. The ideal $(6)$ is $(2,2)$-mpdf but not prime, and $2$ is a zero-divisor modulo $6$. For $n=3$, the primary ideal $(4)$ is $(2,2)$-mpdf and not prime although $3$ is a unit modulo $4$. Example~\ref{ex:xxy} gives an mpdf ideal with prime radical that is not primary. If the radical of an mpdf ideal is maximal, then the ideal is primary, as is every ideal with maximal radical~\cite{AM}.
\end{remark}

\begin{theorem}\label{thm:n2}
Let $I$ be a proper ideal of $R$. The following statements are equivalent.
\begin{enumerate}[label=\textup{(\roman*)}]
\item $I$ is $(2k,2)$-mpdf for some $k\ge1$.
\item $I$ is $(2k,2)$-mpdf for every $k\ge1$.
\item $I$ is radical, and whenever $x,y\in R$ satisfy $xy\in I$, we have $2x\in I$ or $2y\in I$.
\end{enumerate}
\end{theorem}
\begin{proof}
(ii)$\Rightarrow$(i) is trivial, and (i)$\Rightarrow$(iii) is Corollary~\ref{cor:closedradical} together with Lemma~\ref{lem:n2}.

(iii)$\Rightarrow$(ii). Fix $k\ge1$ and let $a,b\in R$ with $a^{2k}-b^2\in I$. We must show $a^k-b\in I$ or $a^k+b\in I$. Put $u=a^k-b$ and $v=a^k+b$. Then $uv=a^{2k}-b^2\in I$ and $u+v=2a^k$. By \textup{(iii)} applied to $uv\in I$, we have $2u\in I$ or $2v\in I$.
\begin{itemize}
\item[Case (1)] If $2u\in I$, then $u^2=u(u+v)-uv=a^k(2u)-uv\in I$, so $u\in I$ because $I$ is radical.
\item[Case (2)] If $2v\in I$, then $v^2=v(u+v)-uv=a^k(2v)-uv\in I$, so $v\in I$.
\end{itemize}
In both cases the $(2k,2)$-mpdf condition holds for $a,b$.
\end{proof}

\begin{corollary}\label{cor:n2}
Let $I$ be a proper ideal of $R$.
\begin{enumerate}
\item If $2\in I$, then $I$ is $(2k,2)$-mpdf if and only if $I$ is radical.
\item If $2$ is a non-zero-divisor on $R/I$, then $I$ is $(2k,2)$-mpdf if and only if $I$ is prime.
\item The property of being $(2k,2)$-mpdf does not depend on $k$.
\end{enumerate}
\end{corollary}
\begin{proof}
\textup{(a)} If $2\in I$, then $2x\in I$ for all $x\in R$, so the second condition in Theorem~\ref{thm:n2}\textup{(iii)} holds automatically. Theorem~\ref{thm:n2} gives the claim.

\noindent\textup{(b)} If $I$ is $(2k,2)$-mpdf, then $I$ is prime by Theorem~\ref{thm:prime}(a). The converse is Proposition~\ref{prop:prime}(a).

\noindent\textup{(c)} This is the equivalence of \textup{(i)} and \textup{(ii)} in Theorem~\ref{thm:n2}.
\end{proof}

The following minimal-prime description specializes, when $k=1$ and
$I\ne0$, to Hizem's characterization~\cite{Hizem}. The preceding
theorem supplies independence of the multiplier under our convention.

\begin{theorem}\label{thm:square-minimal}
A proper ideal $I$ is $(2k,2)$-mpdf if and only if $I$ is radical and $2$ belongs to every prime minimal over $I$, with at most one exception.
\end{theorem}
\begin{proof}
Let $I$ be $(2k,2)$-mpdf. Then $I$ is radical by Corollary~\ref{cor:closedradical}. Let $P\in\Min_R(I)$ with $2\notin P$, and let $Q\in\Min_R(I)$ with $Q\ne P$; we show $2\in Q$.

Since $P$ and $Q$ are minimal over $I$, neither contains the other, so we may choose $x\in P\setminus Q$. The ring $R_P/IR_P$ is reduced because $I$ is radical, and it has a unique prime ideal because $P$ is minimal over $I$; a reduced ring with a unique prime ideal is a field, so $IR_P=PR_P$. In particular $x/1\in IR_P$, and there is $c\in R\setminus P$ with $cx\in I$. By Lemma~\ref{lem:n2}, $2c\in I$ or $2x\in I$. Since $2\notin P$ and $c\notin P$, we have $2c\notin P\supseteq I$, so $2x\in I\subseteq Q$. As $x\notin Q$ and $Q$ is prime, $2\in Q$.

Conversely, suppose $I=\rad I=\bigcap_{P\in\Min_R(I)}P$ and at most one minimal prime $P_0$ does not contain $2$. Let $a,b\in R$ with $a^{2k}-b^2\in I$, and put $u=a^k-b$, $v=a^k+b$, so that $uv\in I$ and $v-u=2b$.
\begin{itemize}
\item[Case (1)] Let $P\in\Min_R(I)$ with $2\in P$. Since $uv\in P$ and $P$ is prime, $u\in P$ or $v\in P$. As $v-u=2b\in P$, both $u$ and $v$ lie in $P$.
\item[Case (2)] If there is no exceptional prime, then $u,v\in P$ for every minimal prime $P$, so $u\in I$.
\item[Case (3)] If $P_0$ exists, then $uv\in P_0$ gives $w\in\{u,v\}$ with $w\in P_0$. By the first item, $w\in P$ for all other minimal primes $P$. So $w\in\bigcap_PP=I$.
\end{itemize}
Thus $u\in I$ or $v=F_{k,2}(a,b)\in I$.
\end{proof}

\begin{remark}
For $k=1$ and nonzero $I$, Proposition~\ref{prop:sdfcdf}(a) and Theorem~\ref{thm:n2} contain the facts from~\cite{ABC} that nonzero sdf-absorbing ideals are radical, that every nonzero radical ideal is sdf-absorbing when $\Char(R)=2$, and that nonzero sdf-absorbing ideals are prime when $2\in R^\times$. Condition (iii) of Theorem~\ref{thm:n2} also covers the case in which $2$ is a nonzero zero-divisor modulo $I$, as for $(6)$ in $\Z$. It says that a product lying in $I$ has a factor in the colon ideal $(I:2)$.
\end{remark}

\begin{proposition}
Let $(R,M)$ be a local ring with $2\notin M$. Then the $(2k,2)$-mpdf ideals of $R$ are exactly the prime ideals of $R$.
\end{proposition}
\begin{proof}
Since $2\notin M$ and $M$ is the unique maximal ideal, $2$ is a unit of $R$, and its image in every nonzero quotient $R/I$ is a unit, hence a non-zero-divisor. Corollary~\ref{cor:n2}(b) applies to every proper ideal.
\end{proof}

\subsection{Even exponents and primary ideals}

The primeness result for the square condition in
\cite[Theorem~2.6]{ABC} has an extension to every even exponent.
The conclusion is that the ideal is primary; nilpotent elements prevent
primeness in general.

\begin{theorem}\label{thm:even-primary}
Let $R$ be a commutative ring, let $I\subsetneq R$ be a
$(kn,n)$-mpdf ideal, where $k\ge1$ and $n\ge2$ is even. If
multiplication by $n$ is injective on $R/I$, then $I$ is primary.
For $n=2$, the ideal $I$ is prime.
\end{theorem}
\begin{proof}
Put $S=R/I$ and $N=\sqrt{0_S}$. We first show that $N$ is prime.
Let $x,y\in S$ with $xy=0$. Set $a=x+y$ and $b=y^k-x^k$.
All mixed terms in the relevant binomial expansions vanish, and
the evenness of $n$ gives
\[
 a^{kn}=b^n,\qquad
 a^k-b=2x^k,\qquad
 F_{k,n}(a,b)=n y^{k(n-1)}.
\]
For the last identity, put $u=x^k$ and $v=y^k$, so $uv=0$;
the pure $u$ terms have total coefficient
$\sum_{i=0}^{n-1}(-1)^i=0$, and the pure $v$ terms have total
coefficient $n$. The mpdf property implies $2x^k=0$ or
$n y^{k(n-1)}=0$. The element $2$ is regular because it divides
the regular element $n$. Hence $x$ or $y$ is nilpotent.
If $xy$ is merely nilpotent, apply this conclusion to $x^h,y^h$
for an integer $h$ with $(xy)^h=0$. It follows that $N$ is prime.

We next show that every zero divisor of $S$ is nilpotent.
Suppose $xy=0$ with $y\ne0$ and $x\notin N$. Since $N$ is
prime, $y\in N$. Choose the least $t\ge2$ with $y^t=0$, and put
$z=y^{t-1}$. Then $z\ne0$, $z^2=0$, and $xz=0$.
For $a=x$ and $b=x^k+z$ we have
\[
 a^{kn}=b^n,\qquad a^k-b=-z\ne0,\qquad
 F_{k,n}(a,b)=n x^{k(n-1)}
     +\binom n2 x^{k(n-2)}z.
\]
The first summand in the cofactor is not nilpotent: if its $h$th
power were zero, regularity of $n^h$ would imply
$x^{hk(n-1)}=0$, contrary to $x\notin N$. The second summand is
nilpotent because it is a multiple of $z$. Their sum is therefore
nonzero. This contradicts the mpdf property, and proves the assertion.
Thus the zero ideal of $S$ is primary, equivalently $I$ is primary.
When $n=2$, Corollary~\ref{cor:closedradical} also makes $I$
radical, so it is prime.
\end{proof}

\begin{remark}\label{rem:regular-prime-radical}
The second part of the proof does not use evenness. It proves, for any
$n\ge2$, that an $(kn,n)$-mpdf ideal $I$ is primary whenever $n$ is a
non-zero-divisor on $R/I$ and $\sqrt I$ is prime. The example of a
nonprimary mpdf ideal with prime radical in
Example~\ref{ex:xxy} has $n=0$ in its quotient.
\end{remark}

\begin{example}\label{ex:even-primary-sharp}
Neither the evenness hypothesis nor the distinction between primary
and prime can be suppressed.
\begin{enumerate}
\item The ideal $10\mathbb Z$ is $(3k,3)$-mpdf for every $k\ge1$
and $3$ is a unit modulo $10$, but $10\mathbb Z$ is not primary.
Indeed, the cube map is injective on both $\mathbb F_2$ and
$\mathbb F_5$, and hence on $\mathbb Z_{10}$ by the Chinese
remainder theorem. However, $2\cdot5\in10\mathbb Z$, while
$2\notin10\mathbb Z$ and no power of $5$ belongs to $10\mathbb Z$.
\item The nonradical primary ideal $(X^2)\subset\mathbb Q[X]$
is $(kn,n)$-mpdf for every $n\ge3$ and $k\ge1$.
To see this, work in $S=\mathbb Q[\epsilon]/(\epsilon^2)$ and
put $x=a^k$, $y=b$ for a pair with $x^n=y^n$. If the residues of
$x,y$ in $\mathbb Q$ differ, $x-y$ is a unit, so $F_n(x,y)=0$.
If their common residue is nonzero, then the residue of $F_n(x,y)$
is $n\overline x^{\,n-1}\ne0$, so this cofactor is a unit and
$x-y=0$. If both residues are zero, then $x,y\in\epsilon S$,
and every monomial of $F_n(x,y)$ has degree $n-1\ge2$, so the
cofactor is zero. Thus being primary cannot be strengthened to
primeness in Theorem~\ref{thm:even-primary} for even $n\ge4$.
Here $F_n(x,y)=\sum_{i=0}^{n-1}x^{n-1-i}y^i$.
\end{enumerate}
\end{example}

\Needspace{4\baselineskip}\subsection{Quartic pair products and primary ideals}
\label{quartic:primary-section}

For the quartic pair-product condition, the nilpotency bound in
\cite[Proposition~4.1]{TAHA} states that every nilpotent has cube zero
when $2$ is a unit. The next lemma strengthens the bound to the cube
of the entire nilradical and requires only regularity of $2$.
Here $\Nil(S)=\sqrt{0_S}$ denotes the nilradical of $S$.

\begin{lemma}\label{quartic:nil-cube}
Let $S$ be a nonzero commutative ring whose zero ideal is
$4$-pdf-absorbing. If multiplication by $2$ is injective on $S$, then
$\Nil(S)^3=0$.
\end{lemma}
\begin{proof}
Put $N=\Nil(S)$. The map $S\to S[1/2]$ is injective, and the zero
ideal of $S[1/2]$ is $4$-pdf-absorbing by
\cite[Proposition~3.6]{TAHA}. Since $2$ is a unit there,
\cite[Proposition~4.1]{TAHA} gives $(x/1)^3=0$ for every $x\in N$.
Injectivity yields $x^3=0$ in $S$.

Let $a,b\in N$ and put $u=a^2b$, $v=ab^2$. Their fourth powers
vanish, and their three pair products are
\[
 a^2-b^2,\qquad v-u,\qquad u+v.
\]
If $a^2=b^2$, then $u=b^3=0$ and $v=a^3=0$. Thus the
$4$-pdf condition gives $u=v=0$, or $v=u$, or $v=-u$.
Apply it also to $a+b$ and $a-b$, whose cubes and fourth powers
vanish because they belong to $N$. Their pair products are
\[
 4ab,\qquad 4b(a^2+b^2)=4u,\qquad
 4a(a^2+b^2)=4v.
\]
Regularity of $2$ gives $ab=0$, or $u=0$, or $v=0$.
Together with the preceding alternatives, this forces $u=v=0$.
Hence $a^2b=0$ for all $a,b\in N$.
Finally, for $x,y,z\in N$,
\[
 2xyz=(x+y)^2z-x^2z-y^2z=0.
\]
Cancellation of $2$ gives $xyz=0$, proving $N^3=0$.
\end{proof}

In particular, if $I$ is $4$-pdf-absorbing and $2$ is regular on
$R/I$, then $(\sqrt I)^3\subseteq I$.
Corollary~4.2 of \cite{TAHA} obtains this containment when $2$ and
$3$ are units modulo $I$, using polarization and cancellation of $6$.
The additional pair $(a+b,a-b)$ above removes the hypothesis on $3$.

The equivalence of \textup{(i)} and \textup{(iii)} in the next theorem
under $2,3\notin P$ is \cite[Theorem~4.6]{TAHA}; its sufficient
direction already assumes only $2\notin P$ in
\cite[Proposition~4.5]{TAHA}. The theorem also identifies the strict
quartic condition under the weaker hypothesis.

\begin{theorem}\label{quartic:primary-equivalence}
Let $I$ be a $P$-primary ideal of $R$ with $2\notin P$.
The following statements are equivalent.
\begin{enumerate}[label=\textup{(\roman*)}]
\item $I$ is $4$-pdf-absorbing.
\item $I$ is $(4,4)$-mpdf.
\item $P^3\subseteq I$.
\end{enumerate}
\end{theorem}
\begin{proof}
Proposition~\ref{quartic:implication} gives (ii)$\Rightarrow$(i).
Since $I$ is $P$-primary and $2\notin P$, multiplication by $2$
is injective on $R/I$. Its nilradical is $P/I$, so
Lemma~\ref{quartic:nil-cube} gives (i)$\Rightarrow$(iii).

Assume (iii), and let $a^4-b^4\in I$. Put
$D=a-b$ and $C=(a+b)(a^2+b^2)$, so $DC\in I$.
If $a,b\in P$, then $C\in P^3\subseteq I$.
Otherwise, if $D\notin P$, primaryness gives $C\in I$.
If $D\in P$, the images of $a,b$ in the domain $R/P$ are equal
and nonzero. Since $2\notin P$, the images of $a+b$ and
$a^2+b^2$ are nonzero. Thus $C\notin P$, and primaryness gives
$D\in I$. This proves (iii)$\Rightarrow$(ii).
\end{proof}

This theorem concerns the unmixed exponent $(4,4)$.
The comparison with $(4k,4)$ for $k\ge2$ is developed using the
valuation criteria in the companion paper~\cite{MPDFII}.

\begin{example}\label{quartic:two-essential}
The restriction $2\notin P$ is essential in both directions.
In $R=\F_2[t]$, put $P=(t)$ and $I=(t^5)$. If
$a^4-b^4=(a-b)^4\in I$, then $a-b\in(t^2)$, so the strict
quartic cofactor $(a-b)^3$ belongs to $(t^6)\subseteq I$.
Thus $I$ is $(4,4)$-mpdf and $4$-pdf-absorbing, but
$P^3\nsubseteq I$.

For the other direction, \cite[Example~4.9]{TAHA} takes
\[
 R=(\Z/4\Z)[t]/(2t^2,t^3),\qquad P=(2,t),\qquad I=0.
\]
Here $P^3=0$ and $R/P\cong\F_2$, so $R$ is local with nilpotent
maximal ideal and $I$ is $P$-primary. For $a=1+t$ and $b=1$,
the fourth powers agree but the pair products are
$2t+t^2$, $2t$, and $2t$, all nonzero. Hence $I$ is not
$4$-pdf-absorbing and therefore not $(4,4)$-mpdf.
\end{example}

\Needspace{4\baselineskip}\subsection{Prime characteristic and \texorpdfstring{$p$}{p}-power exponents}

When the exponent is a power of $p$ and $p\in I$, the Frobenius
identities identify the mpdf condition with a power-closure condition.

\begin{theorem}\label{thm:charp}
Let $p$ be a prime, $e\ge1$ and $n=p^e$, and let $I$ be a proper ideal of $R$ with $p\in I$. The following statements are equivalent.
\begin{enumerate}[label=\textup{(\roman*)}]
\item $I$ is $(kn,n)$-mpdf for some $k\ge 1$.
\item $I$ is $(kn,n)$-mpdf for every $k\ge1$.
\item $I$ is $(n,n-1)$-closed.
\end{enumerate}
In this case $F_{k,n}(a,b)\in I$ whenever $a^{kn}-b^n\in I$.
\end{theorem}
\begin{proof}
(ii)$\Rightarrow$(i) is trivial, and (i)$\Rightarrow$(iii) is Proposition~\ref{prop:special}.

(iii)$\Rightarrow$(ii). Fix $k\ge1$ and let $a,b\in R$ with $a^{kn}-b^n\in I$. Put $c=a^k-b$.

Since $pR\subseteq I$, Lemma~\ref{lem:charp} gives
\[
c^n\equiv a^{kn}-b^n\equiv0\pmod I\qquad\text{and}\qquad F_{k,n}(a,b)\equiv c^{n-1}\pmod I.
\]
Thus $c^n\in I$, and since $I$ is $(n,n-1)$-closed, $c^{n-1}\in I$. Hence $F_{k,n}(a,b)\in I$. So the implication in the definition holds through its second alternative, which also proves the final assertion.
\end{proof}

\begin{remark}
For $n=p=2$ and $n=p=3$, Theorem~\ref{thm:charp} recovers the characteristic-two and characteristic-three criteria in~\cite[Theorem 2.4]{ABC} and~\cite[Theorem 2.6]{Far25}, respectively. Those results assume $\Char(R)=p$; the theorem above only needs $p\in I$.
\end{remark}

\begin{corollary}\label{quartic:characteristic-two}
Let $I$ be a proper ideal of $R$ with $2\in I$. The following
statements are equivalent.
\begin{enumerate}[label=\textup{(\roman*)}]
\item $I$ is $4$-pdf-absorbing.
\item $I$ is $(4,3)$-closed.
\item $I$ is $(4k,4)$-mpdf for some, equivalently every, $k\ge1$.
\end{enumerate}
\end{corollary}
\begin{proof}
In $S=R/I$, put $d=a-b$. The fourth-power difference is $d^4$,
and the three pair products are $d^2,d^3,d^3$.
Thus (i) is equivalent to $d^4=0\Rightarrow d^3=0$ for every
$d\in S$; necessity may be tested on $(a,b)=(d,0)$.
This is the quotient form of \cite[Proposition~2.18]{TAHA}.
The equivalence with (iii) is the case $p=2$, $n=4$ of
Theorem~\ref{thm:charp}.
\end{proof}

For a field $F$ of characteristic two,
\cite[Corollary~2.19]{TAHA} gives the exact list
$e\in\{1,2,3,5,6,9\}$ for which $(t^e)\subset F[t]$ is
$4$-pdf-absorbing. By Corollary~\ref{quartic:characteristic-two},
this is also the list for $(4k,4)$-mpdf absorption, for every $k$.

The prime-power criterion also controls exponents that are not prime powers.

\begin{theorem}\label{thm:divexponent}
Let $p$ be prime, $p\in I$, and $n=p^e$ with $e\ge1$. If $I$ is $(kn,n)$-mpdf for one $k\ge1$, then, for all $q,h\ge1$, it is $(hqn,qn)$-mpdf.
\end{theorem}
\begin{proof}
Let $S=R/I$, a ring in which $p=0$. By Proposition~\ref{prop:special}, $I$ is $(n,n-1)$-closed, so $u^n=0$ implies $u^{n-1}=0$ for every $u\in S$. Let $a,b\in R$ with $a^{hqn}-b^{qn}\in I$, and put $x=\bar a^h$, $y=\bar b$ and $d=x^q-y^q$.

In characteristic $p$ the Frobenius map gives $d^n=(x^q)^n-(y^q)^n=x^{qn}-y^{qn}=0$, because $x^{qn}=\bar a^{hqn}=\bar b^{qn}=y^{qn}$. Hence $d^{n-1}=0$. By Lemma~\ref{lem:poly}(a) and (d), and then Lemma~\ref{lem:charp} in the ring $S$,
\begin{align*}
F_{h,qn}(\bar a,\bar b)=F_{1,qn}(x,y)&=F_{1,q}(x,y)\,F_{1,n}(x^q,y^q)\\
&=F_{1,q}(x,y)\,(x^q-y^q)^{n-1}=F_{1,q}(x,y)\,d^{n-1}=0.
\end{align*}
Thus $F_{h,qn}(a,b)\in I$, and $I$ is $(hqn,qn)$-mpdf.
\end{proof}

The converse along a divisibility chain fails. The ideal $(x^3)\subset\F_3[x]$ is not $(3k,3)$-mpdf for any $k$, since $x^3\in(x^3)$ but $x^2\notin(x^3)$. It is $(9k,9)$-mpdf for every $k$: if $c^9\in(x^3)$, then $x\mid c$, so $c^8\in(x^3)$, and Theorem~\ref{thm:charp} applies.

\begin{proposition}
Let $p\in I$ for a prime $p$, and suppose that $(\rad{I})^t\subseteq I$ for some integer $t\ge1$. Choose $e\ge1$ with $p^e-1\ge t$. Then $I$ is $(hN,N)$-mpdf for every $h\ge1$ and every positive multiple $N$ of $p^e$.
\end{proposition}
\begin{proof}
Let $c\in R$ with $c^{p^e}\in I$. Then $c\in\rad I$, so $c^{p^e-1}\in(\rad I)^{p^e-1}\subseteq(\rad I)^t\subseteq I$, because $p^e-1\ge t$. Thus $I$ is $(p^e,p^e-1)$-closed, and Theorem~\ref{thm:charp} shows that $I$ is $(kp^e,p^e)$-mpdf for every $k$. If $N=qp^e$, Theorem~\ref{thm:divexponent} shows that $I$ is $(hqp^e,qp^e)=(hN,N)$-mpdf for every $h$.
\end{proof}

\begin{example}
The following examples show limitations of power-closure conditions in positive characteristic.
\begin{enumerate}
\item In $\F_p[X]$ no single exponent works for every ideal. Let $n=p^e$. Then $X^n\in(X^n)$ while $X^{n-1}\notin(X^n)$, so $(X^n)$ is not $(n,n-1)$-closed, and Proposition~\ref{prop:special} shows that $(X^n)$ is not $(kn,n)$-mpdf for any $k$.
\item In the non-Noetherian ring $R=\F_p[X_2,X_3,\ldots]$, the ideal $I=(X_j^j:j\ge2)$ is not $(kn,n)$-mpdf for any $n\ge2$ and $k\ge1$. Indeed $X_n^n\in I$, while $X_n^{n-1}\notin I$ because $I$ is a monomial ideal and no generator $X_j^j$ divides $X_n^{n-1}$ by Proposition~\ref{prop:special} .
\end{enumerate}
\end{example}

\section{Homomorphisms, localization and products}\label{sec:hom}

Polynomial evaluation makes contraction and quotient transfer formal.
Localization requires a compatible clearing of denominators for the
two powers. The main issue in direct products is different: the same
factor must vanish in every component. We derive an exact criterion
and its consequence for connectedness of the prime spectrum.

\Needspace{4\baselineskip}\subsection{Quotients, homomorphisms and localization}\label{sec:loc}

For $k=1$ and $n=3$, the quotient and homomorphism transfers below
recover~\cite[Theorem 2.16 and Corollary 2.17]{Far25}. We include
their uniform proofs for arbitrary $n,k$ and for use in the subsequent
constructions. The analogous transfers for $4$-pdf-absorbing ideals
are \cite[Proposition~3.2 and Corollary~3.3]{TAHA}.

\begin{theorem}
Let $J\subseteq I$ be ideals of $R$ with $I$ proper. Then $I$ is an $(m,n)$-mpdf ideal of $R$ if and only if $I/J$ is an $(m,n)$-mpdf ideal of $R/J$.
\end{theorem}
\begin{proof}
The ideal $I/J$ is proper in $R/J$, and $(R/J)/(I/J)\cong R/I$ by the third isomorphism theorem. By Proposition~\ref{prop:reform}, $I$ is $(m.n)$-mpdf in $R$ if and only if the zero ideal of $R/I$ is $(m.n)$-mpdf, and $I/J$ is $(m.n)$-mpdf in $R/J$ if and only if the zero ideal of $(R/J)/(I/J)$ is $(m.n)$-mpdf. The two zero ideals correspond under the isomorphism.
\end{proof}

\begin{theorem}\label{thm:hom}
Let $f\colon R\to T$ be a ring homomorphism.
\begin{enumerate}
\item If $J$ is an $(m,n)$-mpdf ideal of $T$, then $f^{-1}(J)$ is an $(m,n)$-mpdf ideal of $R$.
\item If $f$ is surjective and $I$ is an $(m,n)$-mpdf ideal of $R$ with $\ker f\subseteq I$, then $f(I)$ is an $(m,n)$-mpdf ideal of $T$.
\item If $f$ is surjective, then $J\mapsto f^{-1}(J)$ is a bijection from the set of $(m,n)$-mpdf ideals of $T$ onto the set of $(m,n)$-mpdf ideals of $R$ containing $\ker f$, with inverse $I\mapsto f(I)$.
\end{enumerate}
\end{theorem}
\begin{proof}
\textup{(a)} The set $f^{-1}(J)$ is an ideal of $R$, and it is proper because $f(1)=1\notin J$. Let $a,b\in R$ with $a^m-b^n\in f^{-1}(J)$. Then $f(a)^m-f(b)^n=f(a^m-b^n)\in J$. Since $J$ is $(m.n)$-mpdf, $f(a)^k-f(b)\in J$ or $F_{k,n}(f(a),f(b))\in J$. As $f$ is a ring homomorphism,
\[
f(a^k-b)=f(a)^k-f(b),\qquad f(F_{k,n}(a,b))=F_{k,n}(f(a),f(b)).
\]
Hence $a^k-b\in f^{-1}(J)$ or $F_{k,n}(a,b)\in f^{-1}(J)$.

\noindent\textup{(b)} Since $f$ is surjective, $f(I)$ is an ideal of $T$. It is proper: if $1=f(i)$ with $i\in I$, then $1-i\in\ker f\subseteq I$, so $1\in I$, which is false. The composite $R\to T\to T/f(I)$ is surjective, and its kernel is $f^{-1}(f(I))=I+\ker f=I$. Hence $T/f(I)\cong R/I$, and Proposition~\ref{prop:reform} shows that $f(I)$ is $(m.n)$-mpdf.

\noindent\textup{(c)} For an mpdf ideal $J$ of $T$, $f^{-1}(J)$ is $(m.n)$-mpdf and contains $\ker f$ by \textup{(a)}, and $f(f^{-1}(J))=J$ since $f$ is surjective. For an mpdf ideal $I\supseteq\ker f$, $f(I)$ is $(m.n)$-mpdf by \textup{(b)}, and $f^{-1}(f(I))=I+\ker f=I$. So the two maps are mutually inverse.
\end{proof}

\begin{corollary}\label{cor:sub}
If $R\subseteq T$ is a subring and $J$ is an $(m,n)$-mpdf ideal of $T$, then $J\cap R$ is an $(m,n)$-mpdf ideal of $R$. If $R\to T$ is faithfully flat and $IT$ is an $(m,n)$-mpdf ideal of $T$, then $I$ is an $(m,n)$-mpdf ideal of $R$.
\end{corollary}
\begin{proof}
For the inclusion $\iota\colon R\to T$ we have $\iota^{-1}(J)=J\cap R$, so the first statement is Theorem~\ref{thm:hom}(a). If $R\to T$ is faithfully flat, then $IT\cap R=I$ (see~\cite{AM}); so $I=\iota^{-1}(IT)$ is $(m.n)$-mpdf by Theorem~\ref{thm:hom}(a).
\end{proof}

In particular, if $I[X]$ is $(m.n)$-mpdf in $R[X]$, then $I$ is $(m.n)$-mpdf in $R$. Unlike the corresponding statement for sdf-absorbing ideals in~\cite[Theorem 2.10]{ABC}, Theorem~\ref{thm:hom}(a) needs no hypothesis that the target ideal $J$ be nonzero; this is a consequence of allowing $a=0$ or $b=0$ in Definition~\ref{def:mpdf}.

\begin{example}
Contractions of mpdf ideals are $(m.n)$-mpdf, but extensions need not be, and the hypothesis $\ker f\subseteq I$ in Theorem~\ref{thm:hom}(b) cannot be omitted.
\begin{enumerate}
\item Let $p$ be an odd prime and $f\colon\Z\to\Z\times\Z$ the diagonal embedding. The prime ideal $p\Z$ is $(m.n)$-mpdf, but its extension $p\Z\times p\Z$ is not $(2,2)$-mpdf: for $a=(1,1)$ and $b=(1,-1)$ we have $a^2-b^2=(0,0)\in p\Z\times p\Z$, while $a-b=(0,2)$ and $a+b=(2,0)$ lie outside $p\Z\times p\Z$, because $p\nmid2$.
\item Let $f\colon\Z[x]\to\Z$ be evaluation at $0$. The ideal $(x+4)$ is prime, since $\Z[x]/(x+4)\cong\Z$, hence mpdf. Its image $f((x+4))=4\Z$ is not radical ($2^2\in4\Z$, $2\notin4\Z$), hence not $(2,2)$-mpdf by Corollary~\ref{cor:closedradical}. Here $\ker f=(x)\not\subseteq(x+4)$, because $x\equiv-4\ne0$ modulo $(x+4)$.
\item Extension along the faithfully flat map $\Z\to\Z[X]$ need not preserve the mpdf property either. The ideal $(4)\subset\Z$ is $(6,6)$-mpdf: if $a^6=b^6$ in $\Z_4$, then $a,b$ have the same parity. If both are even, $F_{1,6}(a,b)=0$. If both are odd, either $a=b$, or $b=-a$ and
\[
F_{1,6}(a,b)=a^5\sum_{i=0}^{5}(-1)^i=0.
\]
In $\Z_4[X]$, however, take $a=1$ and $b=1+2X$. Then
\[
a^6=b^6=1,\qquad a-b=-2X\ne0,\qquad
F_{1,6}(a,b)=\sum_{i=0}^{5}(1+2iX)=2+2X\ne0.
\]
Thus $(4)\Z[X]$ is not $(6,6)$-mpdf.
\end{enumerate}
\end{example}

The localization argument extends the cubic result
\cite[Theorem 2.15]{Far25}; the differing exponents require the
weighted scaling $(a,b)\mapsto(ca,c^kb)$. For the quartic
pair-product condition, localization preservation is
\cite[Proposition~3.6]{TAHA}.

\begin{theorem}\label{thm:loc}
Let $I$ be an $(m,n)$-mpdf ideal of $R$ and let $S$ be a multiplicatively closed subset of $R$ with $I\cap S=\emptyset$. Then $S^{-1}I$ is an $(m,n)$-mpdf ideal of $S^{-1}R$.
\end{theorem}
\begin{proof}
Since $I\cap S=\emptyset$, $S^{-1}I$ is a proper ideal of $S^{-1}R$~\cite{AM}. Recall that for $r\in R$ we have $r/1\in S^{-1}I$ if and only if $cr\in I$ for some $c\in S$. Let $A=a/s$ and $B=b/t$ with $a,b\in R$ and $s,t\in S$, and suppose $A^m-B^n\in S^{-1}I$. We must show $A^k-B\in S^{-1}I$ or $F_{k,n}(A,B)\in S^{-1}I$.

Put $w=st\in S$, $a'=at\in R$ and $b'=s^kt^{k-1}b\in R$, so that $a'/1=wA$ and $b'/1=w^kB$. Since $m=kn$,
\[
\frac{a'^m-b'^n}{1}=w^mA^m-w^{kn}B^n=w^m(A^m-B^n)\in S^{-1}I,
\]
and hence $c(a'^m-b'^n)\in I$ for some $c\in S$. Because $(c^k)^n=c^m$,
\[
(ca')^m-(c^kb')^n=c^m(a'^m-b'^n)\in I.
\]
As $I$ is $(m.n)$-mpdf, $(ca')^k-c^kb'\in I$ or $F_{k,n}(ca',c^kb')\in I$. By Lemma~\ref{lem:poly}(b) , these elements equal $c^k(a'^k-b')$ and $c^{k(n-1)}F_{k,n}(a',b')$, respectively. Their images in $S^{-1}R$ are
\[
c^k\bigl((wA)^k-w^kB\bigr)=(cw)^k(A^k-B),\qquad c^{k(n-1)}F_{k,n}(wA,w^kB)=(cw)^{k(n-1)}F_{k,n}(A,B).
\]
Thus $(cw)^k(A^k-B)\in S^{-1}I$ or $(cw)^{k(n-1)}F_{k,n}(A,B)\in S^{-1}I$. The element $cw$ lies in $S$, hence is a unit of $S^{-1}R$, and so $A^k-B\in S^{-1}I$ or $F_{k,n}(A,B)\in S^{-1}I$.
\end{proof}

\begin{corollary}\label{cor:loccontr}
Let $S\subseteq R$ be multiplicatively closed and $J$ an $(m,n)$-mpdf ideal of $S^{-1}R$. Then $\{r\in R: r/1\in J\}$ is an $(m,n)$-mpdf ideal of $R$. In particular, if $I$ is $(m,n)$-mpdf with $I\cap S=\emptyset$, then $\{r\in R: sr\in I\text{ for some }s\in S\}$ is $(m,n)$-mpdf.
\end{corollary}
\begin{proof}
Let $\lambda\colon R\to S^{-1}R$ be the canonical map. Then $\lambda^{-1}(J)=\{r\in R:r/1\in J\}$ is $(m.n)$-mpdf by Theorem~\ref{thm:hom}(a). If $I$ is $(m.n)$-mpdf with $I\cap S=\emptyset$, then $S^{-1}I$ is $(m.n)$-mpdf by Theorem~\ref{thm:loc}, and $\lambda^{-1}(S^{-1}I)=\{r\in R:sr\in I\text{ for some }s\in S\}$.
\end{proof}

\begin{proposition}
Let $I$ be a proper ideal of $R$ and $S$ a multiplicatively closed set with $I\cap S=\emptyset$ such that every element of $S$ is a non-zero-divisor on $R/I$. Then $I$ is $(m,n)$-mpdf if and only if $S^{-1}I$ is $(m,n)$-mpdf.
\end{proposition}
\begin{proof}
If $I$ is $(m.n)$-mpdf, so is $S^{-1}I$ by Theorem~\ref{thm:loc}. Conversely, let $S^{-1}I$ be mpdf. We claim $\lambda^{-1}(S^{-1}I)=I$. The inclusion $I\subseteq\lambda^{-1}(S^{-1}I)$ is clear. If $r\in\lambda^{-1}(S^{-1}I)$, then $sr\in I$ for some $s\in S$, so $(s+I)(r+I)=0$ in $R/I$; as $s+I$ is a non-zero-divisor, $r\in I$. Hence $I=\lambda^{-1}(S^{-1}I)$ is $(m.n)$-mpdf by Corollary~\ref{cor:loccontr}.
\end{proof}

\begin{example}
The mpdf property is not local. The ideal $I=(15)$ of $\Z$ is not $(2,2)$-mpdf (Example~\ref{ex:notmpdf}). Let $q$ be a prime and consider $I\Z_{(q)}$.
\begin{itemize}
\item If $q\ne3,5$, then $15$ is a unit of $\Z_{(q)}$, so $I\Z_{(q)}=\Z_{(q)}$ is not proper.
\item If $q=3$, then $5$ is a unit of $\Z_{(3)}$, so $I\Z_{(3)}=3\Z_{(3)}$; likewise $I\Z_{(5)}=5\Z_{(5)}$.
\end{itemize}
Thus every proper localization of $I$ at a prime is a maximal ideal of a local ring, hence prime, hence mpdf.
\end{example}

\Needspace{4\baselineskip}\subsection{Products, intersections and connected spectra}\label{sec:products}

In a product ring, the hypothesis $a^{kn}-b^n\in I$ holds coordinatewise,
but the conclusion requires one factor to belong to $I$ in every
coordinate. Farshadifar~\cite[Theorem 2.19]{Far25} gives partial
two-factor observations for the cubic absorbing condition.
Ye\c{s}ilot, \c{S}im\c{s}ek and \"Ozel Ay
\cite[Propositions 14 and 15]{YSO} give exact finite-product criteria
for cdf-primary ideals using the same coordinate conflict.
For $4$-pdf-absorbing ideals, \cite[Proposition~3.1]{TAHA}
establishes transfer for $I_1\times R_2$ and $R_1\times I_2$.
Our theorem treats products of proper coordinate ideals for arbitrary
$n,k$ and an index set that may be infinite.

The product criterion becomes more transparent once the pairs that always
select the first factor are characterized. Although the defining mixed
power depends on $k$, this stronger condition does not.

\begin{proposition}\label{prop:linear-injectivity}
Let $S$ be a nonzero commutative ring and let $n\ge2$ and $k\ge1$.
The following conditions are equivalent.
\begin{enumerate}[label=\textup{(\roman*)}]
\item For every $a,b\in S$, the equality $a^{kn}=b^n$ implies $a^k=b$.
\item The map $S\to S$, $x\mapsto x^n$, is injective.
\end{enumerate}
In particular, condition \textup{(i)} is independent of $k$, and either
condition implies that $S$ is reduced.
\end{proposition}
\begin{proof}
Condition \textup{(ii)} immediately gives \textup{(i)} by applying the
injectivity to $(a^k)^n=b^n$.
Assume \textup{(i)}. Taking $a=0$ shows that $z^n=0$ implies $z=0$.
If $z$ is nilpotent, then $z^{n^h}=0$ for some $h\ge1$; repeated use of
this implication gives $z=0$. Thus $S$ is reduced.

Suppose $x^n=y^n$. By \textup{(i)} to $a=x$ and
$b=x^{k-1}y$. Since $(x^{k-1}y)^n=x^{kn}$, we obtain
$x^{k-1}(x-y)=0$. Let $P$ be any prime ideal of $S$. If $x\in P$,
then $y^n=x^n\in P$ gives $y\in P$, and hence $x-y\in P$.
If $x\notin P$, the equality $x^{k-1}(x-y)=0$ also gives
$x-y\in P$. Therefore $x-y$ lies in every prime ideal of $S$,
and it is zero because $S$ is reduced. This proves \textup{(ii)}.
\end{proof}

\begin{theorem}\label{thm:product}
Let $\{R_i\}_{i\in\Lambda}$ be a nonempty family of rings and let $I_i$ be a proper ideal of $R_i$ for each $i$. Then $\prod_iI_i$ is an $(m,n)$-mpdf ideal of $\prod_iR_i$ if and only if the following conditions hold.
\begin{enumerate}[label=\textup{(\roman*)}]
\item Each $I_i$ is an $(m,n)$-mpdf ideal of $R_i$.
\item Whenever $i\in\Lambda$ and there exist $a,b\in R_i$ with $a^m-b^n\in I_i$ but $a^k-b\notin I_i$, one has $n\in I_j$ for every $j\ne i$.
\end{enumerate}
\end{theorem}
\begin{proof}
Write $R=\prod_lR_l$ and $I=\prod_lI_l$. Polynomial expressions in elements of $R$ are computed coordinatewise; in particular $F_{k,n}(0,0)=0$, because $n\ge2$ and $k\ge1$.

Suppose $I$ is $(m.n)$-mpdf.

To prove (i), fix $i$ and let $a_i,b_i\in R_i$ with $a_i^m-b_i^n\in I_i$. Let $a,b\in R$ have coordinates $a_i,b_i$ at $i$ and $0$ elsewhere. Then $a^m-b^n\in I$, so $a^k-b\in I$ or $F_{k,n}(a,b)\in I$, and the coordinate at $i$ gives $a_i^k-b_i\in I_i$ or $F_{k,n}(a_i,b_i)\in I_i$.

To prove (ii), let $i\in\Lambda$ and $a_i,b_i\in R_i$ with $a_i^m-b_i^n\in I_i$ and $a_i^k-b_i\notin I_i$, and suppose, for a contradiction, that $n\notin I_j$ for some $j\ne i$. Let $a,b\in R$ have coordinates $(a_i,b_i)$ at $i$, $(1,1)$ at $j$ and $(0,0)$ elsewhere. Then $a^m-b^n\in I$, since its coordinate at $j$ is $1-1=0$. But $a^k-b\notin I$, because its coordinate at $i$ is $a_i^k-b_i\notin I_i$, and $F_{k,n}(a,b)\notin I$, because its coordinate at $j$ is $F_{k,n}(1,1)=n\notin I_j$. This contradicts the mpdf property of $I$.

Conversely,  assume \textup{(i)} and \textup{(ii)}. Let $a=(a_l)$ and $b=(b_l)$ in $R$ with $a^m-b^n\in I$; then $a_l^m-b_l^n\in I_l$ for every $l$. If $a_l^k-b_l\in I_l$ for all $l$, then $a^k-b\in I$ and we are done. Otherwise fix $i$ with $a_i^k-b_i\notin I_i$. We show $F_{k,n}(a,b)\in I$, that is, $F_{k,n}(a_l,b_l)\in I_l$ for every $l$.

For $l=i$, condition (i) gives $F_{k,n}(a_i,b_i)\in I_i$. Let $l\ne i$. By (ii) applied to $a_i,b_i$, we have $n\in I_l$. If $a_l^k-b_l\notin I_l$, then (i) for $I_l$ gives $F_{k,n}(a_l,b_l)\in I_l$. If $a_l^k-b_l\in I_l$, then $b_l\equiv a_l^k\pmod{I_l}$, and Lemma~\ref{lem:poly}(c) gives
\[
F_{k,n}(a_l,b_l)\equiv F_{k,n}(a_l,a_l^k)=na_l^{k(n-1)}\equiv0\pmod{I_l}.
\]
So $F_{k,n}(a_l,b_l)\in I_l$ for every $l$, that is, $F_{k,n}(a,b)\in I$.
\end{proof}

The finite cubic specialization deserves a precise comparison with
\cite[Propositions 14 and 15]{YSO}. For any cdf-absorbing ideal $J$,
$3\in\sqrt J$ implies $3\in J$. Indeed, in $S=R/J$, suppose $3$
is nonzero and nilpotent, and choose the least $s\ge2$ with $3^s=0$.
Put $u=3^{s-1}\ne0$. Then $u^2=3u=0$, so the pair $a=1$,
$b=1+u$ has equal cubes, while $a-b=-u\ne0$ and
$a^2+ab+b^2=3\ne0$, a contradiction. Therefore, when each factor
is cdf-absorbing, the conditions $3\in\sqrt{I_j}$ in those
cdf-primary product criteria are equivalent to $3\in I_j$.
By Proposition~\ref{prop:linear-injectivity}, their cubing-injectivity
alternative is also exactly the alternative in
Theorem~\ref{thm:product} for $n=3$, $k=1$. Thus on finite families
of cdf-absorbing factors, even nonradical ones, the two compatibility
criteria coincide. The general exponents, mixed multipliers and
arbitrary products are the extensions established here.

The next formulation removes the existential witnesses from
Theorem~\ref{thm:product}. It separates the contribution of nonzero
integer scalars from that of the individual factors.

\begin{corollary}\label{cor:product-explicit}
Let $\{R_i\}_{i\in\Lambda}$ be a nonempty family of nonzero commutative
rings, let $I_i\subsetneq R_i$, and fix $n\ge2$ and $k\ge1$.
Put $S_i=R_i/I_i$ and
$E=\{i\in\Lambda:n\cdot1_{S_i}\ne0\}$.
Then $\prod_i I_i$ is $(kn,n)$-mpdf in $\prod_iR_i$ precisely as follows.
\begin{enumerate}
\item If $E=\varnothing$, every $I_i$ must be $(kn,n)$-mpdf, and this
condition suffices.
\item If $E=\{j\}$, the ideal $I_j$ must be $(kn,n)$-mpdf, and the
map $x\mapsto x^n$ must be injective on $S_i$ for every $i\ne j$;
these conditions suffice.
\item If $E$ contains at least two indices, the map $x\mapsto x^n$
must be injective on every $S_i$, and this condition suffices.
\end{enumerate}
In the last case the criterion is independent of $k$, and every $I_i$
is radical.
\end{corollary}
\begin{proof}
By Proposition~\ref{prop:linear-injectivity}, the witnesses in
Theorem~\ref{thm:product}\textup{(ii)} exist at index $i$ exactly when
the $n$th-power map on $S_i$ is not injective. Thus that condition says
that $E\subseteq\{i\}$ whenever this power map fails to be injective.
If $E=\varnothing$, it imposes no restriction. If $E=\{j\}$, it
requires injectivity at every $i\ne j$. If $E$ has at least two
indices, it requires injectivity at every index. Injectivity implies
the mpdf property by Proposition~\ref{prop:prime}\textup{(b)}, so
the individual mpdf hypotheses that remain are exactly those stated.
Reducedness in the last case follows from
Proposition~\ref{prop:linear-injectivity}.
\end{proof}

\begin{corollary}\label{cor:connected-mpdf}
Let $I$ be a $(kn,n)$-mpdf ideal of a commutative ring $R$, and suppose
that multiplication by $n$ is injective on $S=R/I$. If the map
$x\mapsto x^n$ is not injective on $S$, then $\operatorname{Spec}(S)$
is connected. In particular, if $I$ is not radical, then
$\operatorname{Spec}(R/I)$ is connected.
\end{corollary}
\begin{proof}
If the spectrum were disconnected, there would be an idempotent
$e\in S\setminus\{0,1\}$ and a decomposition
$S\cong eS\times(1-e)S$ into two nonzero rings. The injectivity of
multiplication by $n$ gives $ne\ne0$ and $n(1-e)\ne0$.
Corollary~\ref{cor:product-explicit}\textup{(c)}, applied to the two
zero ideals, would make the $n$th-power map injective on both factors
and hence on $S$, a contradiction. Finally, injectivity of the
$n$th-power map implies reducedness by
Proposition~\ref{prop:linear-injectivity}; a nonradical $I$ therefore
satisfies the stated failure of injectivity.
\end{proof}

\begin{example}\label{ex:nilpotent-product-obstruction}
The zero ideals of $\mathbb Z_4$ and $\mathbb F_5$ are
$(3k,3)$-mpdf for every $k\ge1$, but the zero ideal of
$\mathbb Z_4\times\mathbb F_5$ is not $(3k,3)$-mpdf for any $k$.
Indeed, put $a=(0,1)$ and $b=(2,1)$. Then $a^{3k}=b^3=(0,1)$,
whereas $a^k-b=(2,0)\ne0$ and $F_{k,3}(a,b)=(0,3)\ne0$.
This is the nonreduced obstruction in
Corollary~\ref{cor:product-explicit}: a nonreduced factor is compatible
only with other factors on which $n$ vanishes.
\end{example}

\begin{corollary}\label{cor:factors}
The following statements hold for direct products.
\begin{enumerate}
\item If $\prod_jI_j$ is an $(m,n)$-mpdf ideal of $\prod_jR_j$, then every $I_j$ is $(m,n)$-mpdf.
\item If $I_1$ is a proper ideal of $R_1$, then $I_1\times R_2$ is $(m,n)$-mpdf in $R_1\times R_2$ if and only if $I_1$ is $(m,n)$-mpdf in $R_1$.
\end{enumerate}
\end{corollary}
\begin{proof}
\textup{(a)} is condition \textup{(i)} in the necessity part of Theorem~\ref{thm:product}. For \textup{(b)}, $(R_1\times R_2)/(I_1\times R_2)\cong R_1/I_1$, and Proposition~\ref{prop:reform} applies to both ideals.
\end{proof}

The converse of Corollary~\ref{cor:factors}(a) fails: $(0)$ is prime in $\Z$, but $(0)\times(0)$ is not $(2,2)$-mpdf in $\Z\times\Z$ (Example~\ref{ex:notmpdf}(c)). For $n=2$ the compatibility condition is simple.

\begin{theorem}\label{thm:prod2}
Let $r\ge2$, and let $I_j$ be proper ideals of $R_j$ for $1\le j\le r$. Then $I_1\times\dots\times I_r$ is a $(2k,2)$-mpdf ideal of $R_1\times\dots\times R_r$ if and only if every $I_j$ is $(2k,2)$-mpdf and $2\in I_j$ for all but at most one index $j$.
\end{theorem}
\begin{proof}
Each $I_j$ is $(2k,2)$- mpdf by Corollary~\ref{cor:factors}(a). Suppose $2\notin I_i$ and $2\notin I_j$ for distinct $i,j$. Let $a=(1,\ldots,1)$, and let $b$ have $b_i=-1$ and $b_l=1$ for $l\ne i$. Then $a^{2k}=(1,\ldots,1)=b^2$. The coordinate of $a^k-b$ at $i$ is $1-(-1)=2\notin I_i$, and the coordinate of $a^k+b=F_{k,2}(a,b)$ at $j$ is $1+1=2\notin I_j$. So neither factor lies in the product ideal, a contradiction.

Conversely,  suppose every $I_l$ is $(2k,2)$-mpdf and $2\in I_l$ for all $l\ne j_0$ (for some index $j_0$, which we may choose arbitrarily if $2\in I_l$ for all $l$). Let $a,b$ be elements of the product with $a^{2k}-b^2$ in the product ideal, and put $u=a^k-b$, $v=a^k+b$.
\begin{itemize}
\item Let $l\ne j_0$. Since $u_lv_l\in I_l$ and $I_l$ is $(2k,2)$-mpdf, $u_l\in I_l$ or $v_l\in I_l$. As $v_l-u_l=2b_l\in I_l$, both $u_l$ and $v_l$ lie in $I_l$.
\item In coordinate $j_0$, the mpdf property of $I_{j_0}$ gives $w\in\{u,v\}$ with $w_{j_0}\in I_{j_0}$.
\end{itemize}
By the first item, $w_l\in I_l$ for all $l\ne j_0$, so $w$ lies in the product ideal. Since $F_{k,2}(a,b)=v$, one of $a^k-b$, $F_{k,2}(a,b)$ lies in the product ideal.
\end{proof}

\begin{proposition}
Let $I,J$ be comaximal $(m,n)$-mpdf ideals of $R$, then $I\cap J$ is $(m,n)$-mpdf if and only if condition \textup{(ii)} of Theorem~\ref{thm:product} holds for $R_1=R/I$, $R_2=R/J$ and $I_1=I_2=(0)$. For $n=2$ this means $2\in I$ or $2\in J$.
\end{proposition}
\begin{proof}

Suppose that $I+J=R$, then $R/(I\cap J)\cong R/I\times R/J$ by CRT. By Proposition~\ref{prop:reform}, $I\cap J$ is $(m,n)$-mpdf if and only if $(0)\subset R/I\times R/J$ is. Condition \textup{(i)} of Theorem~\ref{thm:product} holds because $I$ and $J$ are mpdf (Proposition~\ref{prop:reform}), so the criterion is condition \textup{(ii)}. For $n=2$ this is the condition of Theorem~\ref{thm:prod2}, namely that $2=0$ in $R/I$ or in $R/J$.
\end{proof}

\begin{example}
In general $I\cap J$ need not be mpdf: $(3)\cap(5)=(15)\subset\Z$ is not $(2,2)$-mpdf by Example~\ref{ex:notmpdf}(a)
\end{example}

\section{Relations with classical ideals}\label{sec:classical}

We compare $(m,n)$-mpdf ideals with several established classes. An ideal $I$ of aring $R$ is \emph{weakly prime} if $0\ne ab\in I$ implies $a\in I$ or $b\in I$
\cite{AS03}. It is \emph{$2$-absorbing} if for all $a,b,c \in R$ with $abc\in I$ implies that
one of $ab,ac,bc$ belongs to $I$ \cite{Ba07}; more generally, it is
\emph{$r$-absorbing} if, whenever $a_1\cdots a_{r+1}\in I$, a product of $r$
of these elements belongs to $I$ \cite{AB11}. It is
\emph{$1$-absorbing prime} if $xyz\in I$, for nonunits $x,y,z \in R$,
implies $xy\in I$ or $z\in I$ \cite{YNN21}. An ideal is
\emph{semiprimary} if its radical is prime \cite{Fu47},
\emph{strongly irreducible} if $J\cap K\subseteq I$ implies
$J\subseteq I$ or $K\subseteq I$ for ideals $J,K$ \cite{HRR02},
and \emph{primal} if the elements $x$ for which $xy\in I$ for some
$y\notin I$ form an ideal \cite{Fu50}.

Prime ideals are $(m.n)$-mpdf, and mpdf ideals are $(n,n-1)$-closed
(Proposition~\ref{prop:prime} and Corollary~\ref{cor:closedradical}).
The latter implication reverses when $n=p^e$ and $p\in I$
(Theorem~\ref{thm:charp}); a primary square-case mpdf ideal is prime
(Theorem~\ref{thm:prime}). The following comparisons use only the
displayed parameters, without asserting independence for every fixed
$(m,n)$.

\Needspace{20\baselineskip}
\begin{proposition}\label{prop:classical-comparisons}
The examples in the table have the stated properties. Unless another
ring is specified, the ideals belong to $\Z$. For each $r\ge1$, let
$N_r$ be a product of $r+1$ distinct primes congruent to $2$ modulo $3$.

\begin{center}
\normalfont\small
\renewcommand{\arraystretch}{1.12}
\begin{tabularx}{\linewidth}{@{}>{\raggedright\arraybackslash}p{0.24\linewidth}
  >{\raggedright\arraybackslash}X
  >{\raggedright\arraybackslash}X@{}}
\toprule
\textbf{Class} & \textbf{Mpdf, but not in the class}
 & \textbf{In the class, but not mpdf}\\
\midrule
Radical & $(4)$: $(3,3)$ & $(15)$: $(2,2)$\\
Primary & $(6)$: $(2,2)$ & $(9)$: $(3,3)$\\
Semiprimary & $(6)$: $(2,2)$ & $(9)$: $(3,3)$\\
Weakly prime & $(6)$: $(2,2)$ & $(0)\subset\Z\times\Z$: $(2,2)$\\
$2$- and $r$-absorbing
 & $(N_r)$: $(3k,3)$ for all $k\ge1$; not $r$-absorbing
 & $(15)$ is $2$-absorbing: $(2,2)$\\
$1$-absorbing prime
 & $(4)$: $(3,3)$
 & $(0)\subset\Z_9$: $(2,2)$ and $(3,3)$\\
Strongly irreducible; primal
 & $(6)$: $(2,2)$; neither class
 & $(9)$: $(3,3)$; both classes\\
$(n,n-1)$-closed
 & None: mpdf implies this property
 & $(0)\subset\Z\times\Z$: $(2,2)$; $(35)$: $(3,3)$\\
\bottomrule
\end{tabularx}
\end{center}
\end{proposition}
\begin{proof}
The mpdf claims involving $(4),(6),(9),(15)$ and $(0)\subset\Z\times\Z$
were proved in Examples~\ref{ex:six}--\ref{ex:notmpdf}.
The ideals $(15)$ and $(35)$ are radical, whereas $2^2\in(4)$ and
$2\notin(4)$. Also $\rad{(6)}=(6)$ and $\rad{(9)}=(3)$, and
$(9)$ is primary by divisibility by $3$. The product $2\cdot3$,
with no power of $3$ in $(6)$, shows that $(6)$ is neither primary
nor weakly prime. Every zero ideal is weakly prime.

The ideal $(15)$ is $2$-absorbing because the divisors $3$ and $5$
of a product of three elements occur in at most two of its factors.
Primes congruent to $2$ modulo $3$ are infinite: given a finite list,
their product $L$ yields a new such prime divisor of $3L-1$.
For each prime $p\mid N_r$, cubing is bijective on $\F_p$, so
$a^{3k}\equiv b^3\pmod{N_r}$ implies $a^k\equiv b\pmod{N_r}$.
However, the $r+1$ prime factors have product $N_r$, while each
$r$-fold subproduct is $N_r/p\notin(N_r)$.

The nonunits $(x,y,z)=(2,3,2)$ disprove the $1$-absorbing prime
condition for $(4)$. In $\Z_9$, every product of two nonunits is
zero, so its zero ideal is $1$-absorbing prime. It fails the square
mpdf condition because $\bar3$ is nonzero and square-zero, and
the cubic failure follows from that of $(9)\subset\Z$ by
Proposition~\ref{prop:reform}.

The equality $(2)\cap(3)=(6)$ disproves strong irreducibility of
$(6)$. The set in the definition of primality is
$2\Z\cup3\Z$, which is not additive, whereas the corresponding
set for $(9)$ is the ideal $3\Z$. Strong irreducibility of $(9)$
follows from
$v_3(\operatorname{lcm}(a,b))=\max\{v_3(a),v_3(b)\}$ for nonzero
$a,b$; zero generators cause no difficulty.
Finally, radical ideals are $(n,n-1)$-closed. For the cubic failure
of $(35)$, use \cite[Example~2.5]{Far25}, with $a=3, b=-2 $: the mixed difference is $35$,
but the two factors are $5$ and $7$.
\end{proof}

\begin{example}\label{ex:xxy}
In $\F_3[X,Y]$, let $I=(X^2,XY)$. Then $\rad I=(X)$: indeed $X^2\in I$ gives $(X)\subseteq\rad I$, and $I\subseteq(X)$ gives $\rad I\subseteq(X)$. If $c^3\in I$, then $c\in\rad I=(X)$, say $c=Xg$, so $c^2=X^2g^2\in I$. Thus $I$ is $(3,2)$-closed, and since $3=0$ in $\F_3[X,Y]$, Theorem~\ref{thm:charp} shows that $I$ is $(3k,3)$-mpdf for every $k$. Nevertheless $I$ is not primary although its radical is prime: $XY\in I$ and $X\notin I$, but $Y^t\notin I$ for every $t\ge1$, because every element of $I$ is divisible by $X$.
\end{example}

\section{Polynomial and power series extensions}\label{sec:poly}

Contraction gives necessary conditions for adjoining variables. The
converse holds for prime ideals, while the nonradical example in
Section~\ref{sec:hom} shows that a general ascent theorem is impossible.
We distinguish coefficient ideals of formal power series from ideals
generated in the power series ring.
For the quartic pair-product condition, polynomial descent and failure
of ascent from $\Z_8$ to $\Z_8[X]$ are established in
\cite[Proposition~3.9 and Example~3.10]{TAHA}.

For an ideal $I$ of $R$ we write $I[X]=IR[X]$ and $I[[X]]$ for the ideal of power series whose coefficients lie in $I$. Then $R[X]/I[X]\cong(R/I)[X]$ and $R[[X]]/I[[X]]\cong(R/I)[[X]]$, and $I[[X]]=IR[[X]]$ when $I$ is finitely generated.

\begin{proposition}
Let $I\subsetneq R$.
\begin{enumerate}[label=\textup{(\roman*)}]
\item If $I[X]$ or $I[[X]]$ is $(m,n)$-mpdf, then $I$ is $(m,n)$-mpdf.
\item For every $c\in R$, the ideal $(I,X-c)$ of $R[X]$ is $(m,n)$-mpdf if and only if $I$ is.
\end{enumerate}
\end{proposition}
\begin{proof}
\textup{(i)} The contraction of $I[X]$ along the inclusion $R\subseteq R[X]$ is $I$, and likewise for $I[[X]]$ and $R\subseteq R[[X]]$, since an element of $R$ lies in $I[X]$ (respectively $I[[X]]$) if and only if it lies in $I$ (it is a constant polynomial, respectively series, with that coefficient). Corollary~\ref{cor:sub} gives the claim.

\noindent\textup{(ii)} The evaluation $R[X]\to R$, $X\mapsto c$, is surjective, and the composite $R[X]\to R\to R/I$ has kernel $(I,X-c)$. So $R[X]/(I,X-c)\cong R/I$, and Proposition~\ref{prop:reform} applies.

\end{proof}

\section{Pullbacks and related ring constructions}\label{sec:constructions}

This section applies the quotient and product results to familiar ring
constructions. The canonical inverse-image ideals are decided by their
quotients. For intersections of coordinate ideals, compatibility of the
two factor choices is the additional condition.

\Needspace{4\baselineskip}\subsection{Pullback rings}

Let $f:A\to C$ and $g:B\to C$ be surjective homomorphisms. Their pullback is the subring $P=A\times_CB=\{(a,b)\in A\times B:f(a)=g(b)\}$ of $A\times B$, with coordinatewise operations. The coordinate projections $\pi_A:P\to A$ and $\pi_B:P\to B$ and the common-value map $\varepsilon:P\to C$, $\varepsilon(a,b)=f(a)=g(b)$, are surjective. Polynomial identities, in particular Lemma~\ref{lem:poly}(a), hold coordinatewise in $P$; the point requiring care is that one factor must belong to an ideal of $P$ in both coordinates at once. The square-difference case in amalgamations was studied in~\cite{ABC}, and a related weakly $S$-square-difference condition in pullbacks appears in~\cite{AbouhalakaKim}.

\begin{theorem}\label{thm:canonical}
Let $P=A\times_CB$ be as above.
\begin{enumerate}
\item For a proper ideal $I$ of $A$, $\pi_A^{-1}(I)$ is $(m,n)$-mpdf in $P$ if and only if $I$ is $(m,n)$-mpdf in $A$.
\item For a proper ideal $J$ of $B$, $\pi_B^{-1}(J)$ is $(m,n)$-mpdf in $P$ if and only if $J$ is $(m,n)$-mpdf in $B$.
\item For a proper ideal $L$ of $C$, $\varepsilon^{-1}(L)$ is $(m,n)$-mpdf in $P$ if and only if $L$ is $(m,n)$-mpdf in $C$.
\end{enumerate}
\end{theorem}
\begin{proof}
Let $\theta\colon P\to T$ be a surjective homomorphism and $L$ a proper ideal of $T$. The composite $P\to T\to T/L$ is surjective with kernel $\theta^{-1}(L)$, so $\theta^{-1}(L)$ is a proper ideal and $P/\theta^{-1}(L)\cong T/L$. By Proposition~\ref{prop:reform}, $\theta^{-1}(L)$ is $(m.n)$-mpdf in $P$ if and only if the zero ideal of $T/L$ is $(m.n)$-mpdf, if and only if $L$ is $(m.n)$-mpdf in $T$. Apply this to $\theta=\pi_A,\pi_B,\varepsilon$.
\end{proof}

\begin{example}
Let $P=\Z\times_{\Z_6}\Z=\{(a,b):a\equiv b\pmod6\}$ and $K=\ker\varepsilon=6\Z\times6\Z$. Since $P/K\cong\Z_6$ and $(6)$ is $(4,2)$-mpdf (Example~\ref{ex:six}), $K$ is $(4,2)$-mpdf by Theorem~\ref{thm:canonical}(c). However $(2,2),(3,3)\in P$ with $(2,2)(3,3)=(6,6)\in K$ and $(2,2),(3,3)\notin K$, so $K$ is not prime.
\end{example}

For ideals $I\subseteq A$ and $J\subseteq B$, put $H=\{(a,b)\in P:a\in I,\ b\in J\}=P\cap(I\times J)$, an ideal of $P$.

\begin{corollary}
Let $I$ and $J$ be ideals of $A$ and $B$, respectively, with $f(I)=g(J)$, and suppose that $H=P\cap(I\times J)$ is proper.
\begin{enumerate}
\item If $\ker f\subseteq I$, then $H=\pi_B^{-1}(J)$, and $H$ is $(m,n)$-mpdf in $P$ if and only if $J$ is $(m,n)$-mpdf in $B$.
\item If $\ker g\subseteq J$, then $H=\pi_A^{-1}(I)$, and $H$ is $(m,n)$-mpdf in $P$ if and only if $I$ is $(m,n)$-mpdf in $A$.
\end{enumerate}
\end{corollary}
\begin{proof}
\textup{(a)} Assume $\ker f\subseteq I$. Clearly $H\subseteq\pi_B^{-1}(J)$. Conversely, let $(a,b)\in P$ with $b\in J$. Then $f(a)=g(b)\in g(J)=f(I)$, so $f(a)=f(i)$ for some $i\in I$. Hence $a-i\in\ker f\subseteq I$, and $a=(a-i)+i\in I$. Thus $(a,b)\in H$, and $H=\pi_B^{-1}(J)$. The equivalence is Theorem~\ref{thm:canonical}(b). \par\noindent\textup{(b)} Assume $\ker g\subseteq J$. If $(a,b)\in P$ and $a\in I$, then $g(b)=f(a)\in f(I)=g(J)$. Choose $j\in J$ with $g(j)=g(b)$. Then $b-j\in\ker g\subseteq J$, so $b\in J$. Hence $\pi_A^{-1}(I)\subseteq H$; the reverse inclusion follows from the definition of $H$. Thus $H=\pi_A^{-1}(I)$, and Theorem~\ref{thm:canonical}(a) gives the equivalence.
\end{proof}

Even when $I$ and $J$ are $(m.n)$-mpdf, $H$ need not be. The next criterion isolates the only possible conflict between the coordinates.

\begin{proposition}
Assume that $I$ and $J$ are $(m,n)$-mpdf ideals of $A$ and $B$ and that $H$ is proper. Write $u=(a,b)$ and $v=(c,d)$ for elements of $P$. Then $H$ is $(m,n)$-mpdf if and only if no pair $u,v\in P$ with $u^m-v^n\in H$ satisfies one of the following patterns.
\begin{enumerate}
\item $a^k-c\in I$, $F_{k,n}(a,c)\notin I$, $b^k-d\notin J$, $F_{k,n}(b,d)\in J$.
\item $a^k-c\notin I$, $F_{k,n}(a,c)\in I$, $b^k-d\in J$, $F_{k,n}(b,d)\notin J$.
\end{enumerate}
\end{proposition}
\begin{proof}
Let $u=(a,b)$, $v=(c,d)\in P$ with $u^m-v^n\in H$. Then $a^m-c^n\in I$ and $b^m-d^n\in J$.

Since $I$ and $J$ are $(m.n)$-mpdf, $a^k-c\in I$ or $F_{k,n}(a,c)\in I$, and $b^k-d\in J$ or $F_{k,n}(b,d)\in J$. Computing coordinatewise, $u^k-v\in H$ if and only if $a^k-c\in I$ and $b^k-d\in J$, and $F_{k,n}(u,v)\in H$ if and only if $F_{k,n}(a,c)\in I$ and $F_{k,n}(b,d)\in J$.

Suppose both $u^k-v$ and $F_{k,n}(u,v)$ lie outside $H$. If $a^k-c\notin I$, then $F_{k,n}(a,c)\in I$; so $F_{k,n}(b,d)\notin J$ (otherwise $F_{k,n}(u,v)\in H$), and hence $b^k-d\in J$: this is pattern \textup{(b)}. If $a^k-c\in I$, then $b^k-d\notin J$ (otherwise $u^k-v\in H$), so $F_{k,n}(b,d)\in J$, and $F_{k,n}(a,c)\notin I$ (otherwise $F_{k,n}(u,v)\in H$): this is pattern \textup{(a)}. Conversely, each pattern in \textup{(a)}, $b^k-d\notin J$ gives $u^k-v\notin H$ and $F_{k,n}(a,c)\notin I$ gives $F_{k,n}(u,v)\notin H$; in \textup{(b)}, $a^k-c\notin I$ and $F_{k,n}(b,d)\notin J$ do the same.

Hence $H$ fails to be mpdf if and only if some pair $u,v$ with $u^m-v^n\in H$ displays pattern \textup{(a)} or \textup{(b)}.
\end{proof}

\begin{example}
Let $P=\F_5[X]\times_{\F_5}\F_5[Y]$ with constant-term maps. Both coordinate zero ideals are prime. The elements $u=(X,Y)$ and $v=(-X,Y)$ lie in $P$, since all four coordinates have constant term $0$. Then $u^2-v^2=(X^2-X^2,Y^2-Y^2)=0$, but $u-v=(2X,0)\ne0$ and $u+v=(0,2Y)\ne0$. Hence $(0)\subset P$ is not $(2,2)$-mpdf: the coordinates select opposite factors (the first coordinate has linear factor $2X\ne0$ and cofactor $0$;
the second has linear factor $0$ and cofactor $2Y\ne0$). 
\end{example}

\begin{corollary}
If $2$ is a unit in $P$ and $\ker f\ne0\ne\ker g$, then $(0)$ is not $(2k,2)$-mpdf in $P$ for any $k\ge1$.
\end{corollary}
\begin{proof}
Choose $0\ne s\in\ker f$ and $0\ne t\in\ker g$. Then $(s,0)$ and $(0,t)$ lie in $P$, because $f(s)=0=g(0)$ and $f(0)=0=g(t)$. They are nonzero with product $(0,0)=0$, so $(0)$ is not prime. Since $2\in P^\times$ is a non-zero-divisor, Corollary~\ref{cor:n2}(b) shows that $(0)$ is not $(2k,2)$-mpdf.
\end{proof}

\begin{corollary}\label{thm:prime-power}
Let $p$ be a prime dividing $n$, and let $I$ and $J$ be radical ideals of $A$ and $B$ with $p\in I$ and $p\in J$. If $H$ is proper, then $H$ is $(kn,n)$-mpdf for every $k$.
\end{corollary}
\begin{proof}
The kernel of the homomorphism $P\to A/I\times B/J$, $(a,b)\mapsto(a+I,b+J)$, is $H$, so $P/H$ embeds in $A/I\times B/J$. The latter ring is reduced (a product of reduced rings) and $p=0$ in it. Hence $P/H$ is reduced, that is, $H$ is radical, and $p\in H$. Put $S=P/H$ and write $n=pq$. Let $a,b\in S$ with $a^{kn}=b^n$, and set $x=a^k$, $y=b$. Since $S$ is reduced of characteristic $p$,
\[
(x^q-y^q)^p=x^n-y^n=0
\quad\text{gives}\quad x^q=y^q.
\]
By Lemma~\ref{lem:poly}(d),
\[
F_{k,n}(a,b)=F_{1,n}(x,y)
=F_{1,q}(x,y)F_{1,p}(x^q,y^q)
=F_{1,q}(x,y)\,p(x^q)^{p-1}=0.
\]
Proposition~\ref{prop:reform} now gives the conclusion.
\end{proof}

\begin{example}
Let $P=\F_2[X]\times_{\F_2}\F_2[Y]$ with constant-term maps. Taking $I=J=(0)$, Corollary~\ref{thm:prime-power} shows that $(0)\subset P$ is $(2k,2)$-mpdf for every $k$. Yet $(X,0),(0,Y)\in P$ are nonzero with product $0$, so $(0)$ is not prime.
\end{example}

\Needspace{4\baselineskip}\subsection{Idealizations and square-zero extensions}\label{subsec:idealization}

Let $M$ be an $R$-module. Its idealization $R\ltimes M$ has underlying group $R\oplus M$ and multiplication $(a,u)(b,v)=(ab,av+bu)$. For an ideal $I$ of $R$ and a submodule $N$ of $M$, the subgroup $I\oplus N$ is an ideal if and only if $IM\subseteq N$; it is then denoted $I\ltimes N$. These are the homogeneous ideals; not every ideal of $R\ltimes M$ is homogeneous~\cite{AndersonWinders}. When $IM\subseteq N$, there is a natural isomorphism
\begin{equation}\label{eq:idealization-quotient}
(R\ltimes M)/(I\ltimes N)\cong(R/I)\ltimes(M/N).
\end{equation}

The cubic cases of the homogeneous transfers in parts \textup{(i)}
and \textup{(ii)} below occur in
\cite[Theorem 2.20 and Proposition 2.22]{Far25} and
\cite[Theorem 4.6]{Dr26}. The independent quotient proof applies
to every proper $I$, including $I=(0)$. The saturated transfer
$I\ltimes M$ for the quartic pair-product condition is
\cite[Corollary~5.2]{TAHA}.

\begin{theorem}
Let $I\subsetneq R$, $N\subseteq M$ and $IM\subseteq N$.
\begin{enumerate}[label=\textup{(\roman*)}]
\item If $I\ltimes N$ is $(m,n)$-mpdf, then $I$ is $(m,n)$-mpdf.
\item $I\ltimes M$ is $(m,n)$-mpdf if and only if $I$ is.
\item For $n=2$, $I\ltimes N$ is $(2k,2)$-mpdf if and only if $N=M$ and $I$ is $(2k,2)$-mpdf. In fact, every $(2k,2)$-mpdf ideal of $R\ltimes M$ has the form $I\ltimes M$.
\end{enumerate}
\end{theorem}
\begin{proof}
\textup{(i)} The map $\iota\colon R\to R\ltimes M$, $a\mapsto(a,0)$, is a ring homomorphism, and $\iota^{-1}(I\ltimes N)=I$ because $0\in N$. By Theorem~\ref{thm:hom}(a).

\noindent\textup{(ii)} The map $(a,u)\mapsto a+I$ is a surjective homomorphism $R\ltimes M\to R/I$ with kernel $I\ltimes M$. So $(R\ltimes M)/(I\ltimes M)\cong R/I$, and Proposition~\ref{prop:reform} applies.

\noindent\textup{(iii)} Let $L$ be any $(2k,2)$-mpdf ideal of $R\ltimes M$.

$L$ is radical by Corollary~\ref{cor:closedradical}. For $u\in M$ we have $(0,u)^2=(0,0)\in L$, so $(0,u)\in L$; thus $0\ltimes M\subseteq L$. Let $I_0$ be the image of $L$ under the projection $\pi\colon R\ltimes M\to R$, an ideal of $R$. It is proper: if $1\in I_0$, then $(1,u)\in L$ for some $u$, and $(1,u)$ is a unit with inverse $(1,-u)$, contradicting $L\ne R\ltimes M$. Since $L\supseteq\ker\pi=0\ltimes M$, we get $L=\pi^{-1}(I_0)=I_0\ltimes M$, and $I_0$ is $(2k,2)$-mpdf by \textup{(ii)}. Consequently, if $I\ltimes N$ is $(2k,2)$-mpdf, then $I\ltimes N=I_0\ltimes M$, so $N=M$ and $I=I_0$ is $(2k,2)$-mpdf. Conversely, if $N=M$ and $I$ is $(2k,2)$-mpdf, then $I\ltimes M$ is $(2k,2)$-mpdf by \textup{(ii)}.
\end{proof}

The result for $n=2$ parallels~\cite[Theorem 4.16]{ABC}; the all-elements convention also covers $I=(0)$. For $n\ge3$, the equality $N=M$ is not necessary. For example, in
$\mathbb F_3\ltimes M$ every element whose cube is zero has square zero,
so the zero ideal is $(3k,3)$-mpdf for every $k$ by
Theorem~\ref{thm:charp}. A classification of general square-zero
extensions requires additional module and extension data.
For the $4$-pdf condition, \cite[Theorem~5.5]{TAHA} treats
homogeneous ideals with a primary base ideal, $2,3\notin\sqrt I$,
and regular action of $R\setminus\sqrt I$ on the module quotient; the base ideal need not
be prime. The companion paper~\cite{MPDFII} develops the homogeneous
quartic comparison using Theorem~\ref{quartic:primary-equivalence},
and classifies square-zero extensions and nonhomogeneous idealization
ideals with prime first-coordinate projection for $n\ge3$ and $k\ge1$.

\Needspace{4\baselineskip}\subsection{Amalgamated algebras and duplications}\label{subsec:amalgamation}

Let $f:A\to B$ be a homomorphism and $J$ an ideal of $B$. Following~\cite{DAnna}, put $C=A\bowtie^fJ=\{(a,f(a)+j):a\in A,\ j\in J\}\subseteq A\times B$. The coordinate maps $p_1:C\to A$ and $p_2:C\to f(A)+J$ are surjective ring homomorphisms; the second is considered with values in the subring $f(A)+J$ of $B$.

\begin{theorem}
Let $I\subsetneq A$ and $H\subsetneq f(A)+J$ be ideals.
\begin{enumerate}[label=\textup{(\roman*)}]
\item The ideal $I\bowtie^fJ=\{(i,f(i)+j):i\in I,\ j\in J\}$ is $(m,n)$-mpdf in $C$ if and only if $I$ is $(m,n)$-mpdf in $A$.
\item The ideal $p_2^{-1}(H)$ is $(m,n)$-mpdf in $C$ if and only if $H$ is $(m,n)$-mpdf in $f(A)+J$.
\item Suppose $I$ and $H$ are $(m.n)$-mpdf. Then $L=(I\bowtie^fJ)\cap p_2^{-1}(H)$ is $(m.n)$-mpdf provided both of the following hold.
\begin{enumerate}[label=\textup{(\arabic*)}]
\item If there are $a,b\in A$ with $a^m-b^n\in I$ and $a^k-b\notin I$, then $n\in H$.
\item If there are $a,b\in f(A)+J$ with $a^m-b^n\in H$ and $a^k-b\notin H$, then $n\in I$.
\end{enumerate}
In particular, it suffices that $n\in I$ and $n\in H$.
\end{enumerate}
\end{theorem}
\begin{proof}
\textup{(i), (ii)} We have $I\bowtie^fJ=p_1^{-1}(I)$, so $C/(I\bowtie^fJ)\cong A/I$, and $C/p_2^{-1}(H)\cong(f(A)+J)/H$ (Theorem~\ref{thm:canonical} for surjections). Proposition~\ref{prop:reform} gives both equivalences.

\noindent\textup{(iii)} The conditions \textup{(1)} and \textup{(2)} are condition \textup{(ii)} of Theorem~\ref{thm:product} for the ideal $I\times H$ of $A\times(f(A)+J)$ (with $R_1=A$, $R_2=f(A)+J$), and condition \textup{(i)} holds because $I$ and $H$ are $(m.n)$-mpdf. So $I\times H$ is $(m.n)$-mpdf. The ring $C$ is a subring of $A\times(f(A)+J)$, and the contraction of $I\times H$ to $C$ is $\{(a,f(a)+j):a\in I,\ f(a)+j\in H\}=L$. Corollary~\ref{cor:sub} shows that $L$ is $(m.n)$-mpdf. If $n\in I$ and $n\in H$, then \textup{(1)} and \textup{(2)} hold trivially.
\end{proof}

For $k=1$ and $n=3$, part \textup{(i)} recovers the canonical
amalgamation transfer of~\cite[Theorem 2.23]{Far25}; the quotient
proof also covers $I=(0)$. The square and prime-exponent antecedents
are~\cite[Theorem 4.19]{ABC} and~\cite[Theorem 6]{Yes26}.
Part \textup{(iii)} is a sufficient criterion because the two
coordinates of $C$ are linked through $f$ and $J$.

\end{document}